\documentclass[a4paper,reqno,10pt]{amsart}

\usepackage{latexsym}
\usepackage{amssymb}
\usepackage{amscd}
\usepackage{amsmath}
\usepackage{amsfonts}
\usepackage{amsthm}
\usepackage{color}
\usepackage{epsfig}
\usepackage{graphicx}
\usepackage[all]{xy}

\theoremstyle{plain}
\newtheorem{thm}{Theorem}[section]
\newtheorem{cor}[thm]{Corollary}
\newtheorem{lem}[thm]{Lemma}
\newtheorem{prop}[thm]{Proposition}

\newtheorem{ques}[thm]{Question}

\theoremstyle{definition}
\newtheorem{defn}[thm]{Definition}

\newtheorem{rem}[thm]{Remark}

\newtheorem{ex}[thm]{Example}

\newcommand{\C}{{\mathbb C}}
\newcommand{\R}{{\mathbb R}}
\newcommand{\Q}{{\mathbb Q}}
\newcommand{\Z}{{\mathbb Z}}

\def\t{\mathfrak t}

\def\CL{\mathcal{L}}

\def\CN{\mathcal{N}}
\def\CO{\mathcal{O}}

\def\CS{\mathcal{S}}

\def\ep{\varepsilon}
\def\om{\omega}

\def\lb{\lambda}

\def\<{\left\langle}
\def\>{\right\rangle}

\def\i{\sqrt{-1}}

\def\Aut{\operatorname{Aut}}
\def\Lie{\operatorname{Lie}}

\def\Im{\operatorname{Im}}
\def\id{\operatorname{id}}

\def\GL{\operatorname{GL}}
\def\pr{\operatorname{pr}}

\newcommand{\abs}[1]{\lvert#1\rvert}

\numberwithin{equation}{section}

\title{twisted toric structures}
\author{Takahiko Yoshida} 
\thanks{
\ \ \ The author is supported by Research Fellowship of the Japan Society for the 
Promotion of Science for Young Scientists.}

\begin{document}
\begin{abstract}
This paper introduces the notion of twisted toric manifolds which is a generalization of one of 
symplectic toric manifolds, and proves the weak Delzant type classification theorem for them. 
The computation methods for their fundamental groups, cohomology groups in general cases, 
and signatures in four-dimensional cases are also given. 
\end{abstract}
\maketitle
\tableofcontents

\section{Introduction}
 By Delzant's classification theorem \cite{D}, there is a one-to-one correspondence 
between a symplectic toric manifold which is one of the special objects in the theory of 
Hamiltonian torus actions and a Delzant polytope which is a combinatorial object. 
Through this correspondence, various researches on the relationship between symplectic 
geometry, topology, and transformation groups with combinatorics have been done 
\cite{Au, D, G}. 

On the other hand, there exists a manifold such that it may not be itself a symplectic 
toric manifold, but it has a toric structure in a neighborhood of each point all of which 
are patched together in some weak condition. In this paper, as a formulation of such 
manifolds, we shall introduce the notion of twisted toric manifolds and generalize the 
weak version of Delzant's 
classification theorem to them. Recently, some generalizations are also considered 
\cite{HM, K, NS, S, V}. We also investigate the topology of twisted toric manifolds. 
As a result, we can see that there are examples of twisted toric manifolds 
which are not complete non-singular toric varieties in the original algebro-geometric 
sense. In particular, these are not symplectic toric manifolds. 

In general, a twisted toric manifold no longer has a global torus action like that of a 
original symplectic toric manifold, but it has a torus action on a neighborhood of each 
point which comes from a local toric structure and they are patched together in 
certain sense. One of our motivation is to generalize the topological theory of 
transformation group to such a {\it twisted} torus action. Some invariants for transformation 
groups such as equivariant cohomology groups can be generalized to this case and 
we are investigating their properties, in particular, localizations. Unfortunately 
we could not describe this topic in this paper. This will appear later on. 

This paper is organized as follows. First, we recall Hamiltonian torus actions in Section 2 
and symplectic toric manifolds in Section 3 in order that the paper is self-contained. 
Then we shall give the definition of the twisted toric manifold and some examples in 
Section 4. Section 5 is devoted to the classification of twisted toric manifolds. 
In Section 6, We shall compute their fundamental groups and cohomology groups. 
We shall also compute their signatures in four-dimensional cases. 

In the rest of this paper, we shall assume that all manifolds are compact, connected, 
and oriented and all maps preserve orientations, unless otherwise stated. 
\vspace*{3mm}\\
{\bf Acknowledgment:} The author is thankful to Professor Mikio Furuta, 
who is my mentor, for useful comments and suggestions. This work is motivated by 
his suggestion. The author is also thankful to Professor Hisaaki Endo 
for teaching me Meyer's signature cocycle. 

\section{Hamiltonian torus actions}
\subsection{Moment maps}
A {\it symplectic manifold} $(X, \om)$ is a smooth manifold $X$ equipped with 
a non-degenerate closed 2-form $\om$. 
Let us assume that a $k$-dimensional torus $T^k$ acts on $X$ which preserves $\om$. 
In this paper, we identify $T^k$ with $\R^k /\Z^k$, and its Lie algebra $\t$ with $\R^k$. 
By the natural inner product $\<, \>$ on $\R^k$, we also identify the dual space $\t^*$ 
of $\t$ with $\t$ itself. 
\begin{defn}[\cite{Au, GS2}]\label{1.1}
A {\it moment map} for the $T^k$-action is a map $\mu :X\to \t^*$ which is $T^k$-invariant 
with respect to the given $T^k$-action on $X$ and satisfies the condition 
\[
\iota (v_{\xi})\om =\< d \mu , \xi  \>
\]
for $\xi \in \t$, where $v_{\xi}$ is the infinitesimal action, that is, 
the vector field which is defined by 
\begin{equation*}
v_{\xi}(x)=\frac{d}{d\tau }\Big|_{\tau =0}e^{2\pi \i \tau \xi}\cdot x .
\end{equation*}
Note that a moment map for a $T^k$-action is determined up to an additive constant. 
The torus action which has a moment map is said to be \it{Hamiltonian}. 
\end{defn}
Although these are not compact, the following examples are fundamental in this paper.  
\begin{ex}\label{1.2}
Let $(\C^n, \om_{\C^n})$ be the $n$-dimensional complex vector space with 
the symplectic form 
$\displaystyle \om_{\C^n}=\dfrac{-\i }{2\pi }\sum_{i=1}^ndz_i\wedge d\bar{z}_i$. 
$T^n$ acts on $\C^n$ by 
\[
t \cdot z=(e^{2\pi \i t_1}z_1,\ldots ,e^{2\pi \i t_n}z_n) 
\]
for $t =(t_1,\ldots ,t_n)\in T^n$ and $z=(z_1,\ldots ,z_n)\in \C^n$. This action is Hamiltonian 
and a moment map $\mu_{\C^n}:\C^n \to \t^*$ is defined by 
\[
\mu_{\C^n}(z)=(\abs{z_1}^2, \ldots , \abs{z_n}^2). 
\]
In particular, the image of $\mu_{\C^n}$ is 
\[
\t^*_{\ge 0}=\{ \xi =(\xi_1, \ldots , \xi_n)\in \t^*=\R^n \colon 
\xi_i\ge 0\ \text{for}\ i=1, \ldots , n\} .
\] 
\end{ex}
\begin{ex}\label{1.3}
Let $T^*T^n$ be the cotangent bundle of $T^n$. On $T^*T^n$, we fix the natural 
trivialization $T^*T^n\cong \t^* \times T^n$ and identify $T^*T^n$ with 
$\t^*\times T^n$ by this trivialization. $T^*T^n$ has a symplectic form 
$\displaystyle \om_{T^*T^n}=\sum_{i=1}^nd\theta_i \wedge d\xi_i$, where 
$\theta =(\theta_1, \ldots ,\theta_n)$ and $\xi =(\xi_1,\ldots ,\xi_n)$ denote 
the standard coordinates of $T^n$ and $\t^*\cong \R^n$, respectively. $T^n$ acts on 
$T^*T^n$ by 
\[
t \cdot (\xi, \theta )=(\xi, \theta +t )
\] 
for $t \in T^n$ and $(\xi, \theta )\in T^*T^n$. This action is Hamiltonian 
and a moment map $\mu_{T^*T^n}:T^* T^n \to \t^*$ is defined by 
\[
\mu_{T^*T^n}(\xi, \theta )=\xi . 
\]
\end{ex}

\subsection{Symplectic reduction}
Let $(X, \om )$ be a $2n$-dimensional symplectic manifold equipped with a Hamiltonian 
$T^k$-action with a moment map $\mu :X\to \t^*$. There is a method, so called a 
{\it symplectic reduction}, to construct a new symplectic manifold which we shall explain. 
See \cite{Au, GS2} for more details. 
\begin{prop}\label{1.6}
Let $T^k_x (\subset T^k)$ be the stabilizer of $x\in X$. Then the 
annihilator $(\Im d\mu_x)^{\bot}$ of the image of $d\mu_x :T_xX\to \t^*$ 
is isomorphic to the Lie algebra $\t_x$ of $T^k_x$.   
\end{prop}
\begin{proof}
It is clear from the condition in Definition \ref{1.1} and the 
non-degeneracy of $\om$. 
\end{proof}

Let $\varepsilon \in \t^*$. 
Since $\mu$ is invariant under the action, $T^k$-action preserves 
$\mu^{-1}(\varepsilon )$. Suppose that $T^k$-action on $\mu^{-1}(\varepsilon )$ is free. 
Then Proposition \ref{1.6} implies that the level set $\mu^{-1}(\varepsilon )$ 
is smooth, and the quotient space $\mu^{-1}(\varepsilon )/T^k$ is a $2(n-k)$-dimensional 
smooth manifold. In this case, the following proposition is well known.   
\begin{prop}[\cite{MW}]\label{1.7}
The quotient space $\mu^{-1}(\varepsilon )/T^k$ carries a natural symplectic form 
$\om_{\varepsilon}$ such that the equality $\iota^*\om =\pi^* \om_{\varepsilon}$ holds, 
where $\iota$ is a natural inclusion and $\pi$ is a projection 
\[
\xymatrix{
(\mu^{-1}(\varepsilon ), \iota^*\om )\ar[d]_{\pi} \ar@{^{(}->}[r]^{\iota} & (X, \om ) \\
(\mu^{-1}(\varepsilon )/T^k, \om_{\varepsilon}). &   
}
\]
\end{prop}
$(\mu^{-1}(\varepsilon )/T^k, \om_{\varepsilon })$ is called a {\it symplectic quotient}. 

\subsection{Symplectic cutting}
Let us recall the symplectic cutting by Lerman \cite{L}. Suppose that $(X, \om )$ is a 
symplectic manifold equipped with a Hamiltonian $S^1$-action with a moment map 
$\mu :X\to \R$. Define the $S^1$-action on the product space 
$(X\times \C , \om \oplus \om_{\C})$ by 
\[
t\cdot (x, z) =(t\cdot x, e^{-2\pi \sqrt{-1}t}z)
\] 
for $t\in S^1$ and $(x, z)\in X\times \C$. This action is Hamiltonian and the moment map 
$\Phi :X\times \C \to \R$ is  
\[
\Phi (x, z)=\mu (x)-\mu_{\C}(z)=\mu (x)-\abs{z}^2 .
\]  
\begin{prop}\label{1.8}
Let $\varepsilon \in \R$. the $S^1$-action on $\Phi^{-1}(\varepsilon )$ is free, if 
and only if the $S^1$-action on $\mu^{-1}(\varepsilon )$ is free. 
\end{prop}
\begin{proof}
Let $(x,z) \in \Phi^{-1}(\ep )$. If $z\neq 0$, then the stabilizer of $(x, z)$ for 
the $S^1$-action on $X\times \C$ only consists of the unit element since the 
stabilizer of $z$ for the $S^1$-action on $\C$ only consists of the unit 
element. In the case where $z=0$, the stabilizer of $(x, z)$ for the $S^1$-action 
on $X\times \C$ is equal to that of $x$ for the $S^1$-action on $X$. 
This proves the proposition. 
\end{proof}

Assume that the $S^1$-action on $\Phi^{-1}(\varepsilon )$ is free. Then the reduced space 
$\Phi^{-1}(\ep )/S^1$ is a smooth manifold whose dimension is equal to that of $X$. 
Let us consider the reduced space $\Phi^{-1}(\ep )/S^1$. 
The level set $\Phi^{-1}(\ep )$ is a disjoint union of two $S^1$ invariant parts 
\[
\Phi^{-1}(\ep )=\left\{ (x,z)\in X\times \C \colon \mu (x)>\ep ,\ \ 
\abs{z}^2=\mu (x)-\ep \right\} \amalg \mu^{-1}(\ep )\times \{ 0\} . 
\]  
The first part is equivariantly diffeomorphic to the product 
$\{ x\in X \colon$ $ \mu (x)>\ep  \} \times S^1$, and the second part is 
naturally identified with $\mu^{-1}(\ep )$. Then as a set, the quotient space 
$\Phi^{-1}(\ep )/S^1$ is the disjoint union  
\[
\Phi^{-1}(\ep )/S^1 \cong 
\left\{ x\in X \colon \mu (x)>\ep \right\} \amalg \mu^{-1}(\ep )/S^1. 
\]
We denote $\Phi^{-1}(\ep )/S^1$ by $\displaystyle \overline{X}_{\mu \ge \ep }$. 
Topologically, $\displaystyle \overline{X}_{\mu \ge \ep }$ is the quotient of the manifold 
$\displaystyle X_{\mu \ge \ep }=\left\{ x\in X \colon \mu (x)\ge \ep \right\}$ with 
the boundary $\mu^{-1}(\ep )$ by the relation $\sim$, where 
$x\sim x'$ if and only if $x, x'\in \mu^{-1}(\ep )$ and $x'=t\cdot x$ for some 
$t\in S^1$. For this reason, we would like to call the operation that produces 
$\overline{X}_{\mu \ge \ep }$ from the Hamiltonian $S^1$-action on $(X,\om )$ 
{\it symplectic cutting} and $\overline{X}_{\mu \ge \ep }$ is called a {\it cut space}. 
\begin{rem}\label{1.10}
Suppose that $(X,\om )$ has another Hamiltonian $T^k$-action which commutes with 
the $S^1$-action. Then the $T^k$-action on $X$ induces the Hamiltonian $T^k$-action on 
$\overline{X}_{\mu \ge 0}$ simply by letting it act on the first factor. 
\end{rem}
\begin{ex}\label{1.11}
Let us consider Example \ref{1.3} for $n=1$. By Remark \ref{1.10}, the cut space 
$\overline{T^*S^1}_{\mu_{T^*S^1} \ge 0}$ has the Hamiltonian circle action 
which is induced by the original circle action on $T^*S^1$. In this case, 
$\overline{T^*S^1}_{\mu_{T^*S^1} \ge 0}$ is equivariantly symplectomorphic to 
$(\C ,\om_{\C})$ with the circle action in Example \ref{1.2} 
by the symplectomorphism 
$\varphi :\C \to \overline{T^*S^1}_{\mu_{T^*S^1} \ge 0}$ which is defined by 
\[
\varphi (z)=[\abs{z}^2, \dfrac{\arg z}{2\pi}, \abs{z}]. 
\]
\end{ex}
\begin{ex}\label{1.12}
Let $u\in \Z^n\subset \t$, and consider the Hamiltonian $S^1$-action on 
$(T^*T^n, \om_{T^*T^n})$ which is defined by 
\[
t\cdot (\xi, \theta )=(\xi, \theta +tu)
\]
for $t\in S^1=\R /\Z$ and $(\xi, \theta )\in T^*T^n$. The moment map 
$\mu_u :T^*T^n \to \R$ is obtained by 
\[
\mu_u (\xi, \theta )=\< u, \xi \> .
\]
In this case, 
\[
\overline{(T^*T^n)}_{\mu_u \ge \ep }
\cong \left\{ \xi \in \t^* \colon \< u, \xi \> >\ep \right\} 
\amalg \left\{ \xi \in \t^* \colon \< u, \xi \> =\ep \right\} \times T^n/S^1_u , 
\]
where $S^1_u$ is the circle subgroup of $T^n$ generated by $u$. Note that 
$\overline{(T^*T^n)}_{\mu_u \ge \ep }$ is smooth, if and only if $u$ is primitive 
in the sense of Definition \ref{2.2}. For more details, see Appendix A.  
In this case, since the Hamiltonian $T^n$-action 
in example \ref{1.3} commutes with the $S^1$-action, by Remark \ref{1.10}, 
the Hamiltonian $T^n$-action is induced to $\overline{(T^*T^n)}_{\mu_u \ge \ep }$ 
with the moment map 
\[
\overline{\mu_{T^*T^n}}([\xi, \theta , z])=\xi .
\] 
Figure \ref{fig1} shows the change of the image of the moment map under 
symplectic cutting for $n=2$. 
\begin{figure}[hbtp]
\begin{center}
\input{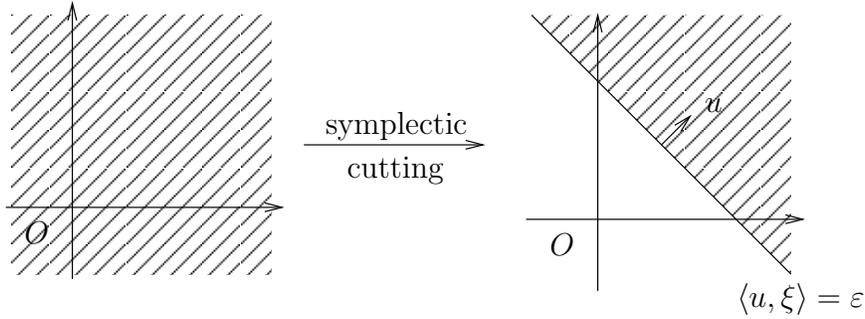}
\caption{the change of the moment image by symplectic cutting}
\label{fig1}
\end{center}
\end{figure}
\end{ex}
\begin{rem}[Simultaneous symplectic cuttings]\label{1.9}
Suppose that $(X, \om )$ is equipped with two commutative Hamiltonian 
$S^1$-actions on $(X,\om )$ with moment maps $\mu_1$ and $\mu_2$.   
Then corresponding symplectic cutting operations also commute each other, and the both 
cut spaces $\overline{(\overline{X}_{\mu_1\ge \ep_1})}_{\mu_2\ge \ep_2}$ 
and $\overline{(\overline{X}_{\mu_2\ge \ep_2})}_{\mu_1\ge \ep_1}$ 
are naturally symplectomorphic to the cut space 
\[
\overline{X}_{\{ \mu_i\ge \ep_i\}_{i=1, 2}}= 
\{ (x,z)\in X\times \C^2 \colon \mu_i(x)-\abs{z_i}^2=\ep_i ,\ i=1, 2\} /T^2 
\]
of the simultaneous symplectic cuttings, that is, the symplectic quotient of the Hamiltonian 
$T^2$-action on $(X\times \C^2, \om \oplus \om_{\C^2})$ with the moment map 
\[
\Phi (x, z)=(\mu_1(x)-\abs{z_1}^2, \mu_2(x)-\abs{z_2}^2)
\]
which is obtained by putting two $S^1$-actions together. 

More generally, for the case where $(X, \om )$ is equipped with $k$ commutative Hamiltonian 
$S^1$-actions with moment maps $\mu_i$ for $i=1, \ldots , k$, the argument goes 
similar way, and the cut spaces are naturally symplectomorphic to the simultaneous cut space 
\[
\overline{X}_{\{ \mu_i\ge \ep_i\}_{i=1,\ldots , k}}= 
\{ (x,z)\in X\times \C^k \colon \mu_i(x)-\abs{z_i}^2=\ep_i ,\ i=1,\ldots ,k \} /T^k.  
\] 
\end{rem}

\section{Symplectic toric manifolds}
For Hamiltonian torus actions, the following fact is well known. 
\begin{thm}[\cite{G}]\label{1.4}
If a $k$-dimensional torus $T^k$ acts effectively on a $2n$-dimensional symplectic 
manifold $(X, \om )$ in a Hamiltonian fashion, then $k\le n$. 
\end{thm}

In particular, in the maximal case of Theorem \ref{1.4}, that is, a closed, connected 
$2n$-dimensional symplectic manifold $(X, \om )$ equipped with an effective 
Hamiltonian $T^n$-action is called a {\it symplectic toric manifold}. 

For a general Hamiltonian torus action, Atiyah and Guillemin-Sternberg show the following 
convexity theorem. 
\begin{thm}[\cite{At, GS1}]
Let $(X, \om )$ be a closed, connected $2n$-dimensional symplectic manifold equipped 
with a Hamiltonian action of $k$-dimensional torus. (We do not require $k=n$.)   
Then the image of a moment map is a convex hull of images of fixed points. 
\end{thm}

In the theory of symplectic toric manifolds, the image of a moment map plays a crucial 
role. 
\begin{defn}\label{2.2}
Let $\{ u_1, \ldots , u_k\}$ is a tuple of vectors of $\Z^n$. $\{ u_1, \ldots , u_k\}$ 
is said to be {\it primitive}, if the sub-lattice ${\rm span}_{\Z}\{ u_1, \ldots , u_k\}$ 
spanned by $u_1, \ldots , u_k$ is a rank $k$ direct summand of the free $\Z$-module 
$\Z^n$. We also say that the tuple $\{ L_1,\ldots , L_k \}$ of rank one sub-lattices in $\Z^n$ 
is {\it primitive}, if there exists a primitive tuple $\{ u_1, \ldots ,u_k\}$ of vectors in $\Z^n$ 
such that each $L_i$ is spanned by $u_i$ for $i=1, \ldots ,k$.   
\end{defn}
\begin{rem}\label{2.25}
The notion of the primitivity of a tuple $\{ u_1, \ldots , u_k \}$ of vectors (hence the 
tuple $\{ \Z_1,\ldots , \Z_k\}$ of rank one sub-lattices) in $\Z^n$ is invariant under the 
action of $GL_n(\Z )$. 
\end{rem}

Let $\Delta$ be a convex polytope in $\t^*(=\R^n)$ in $\t^*$ which is written by 
\begin{equation}\label{eq2.1}
\Delta =\bigcap_{i=1}^d\left\{  \xi \in \t^* \colon \< u_i, \xi \> \ge \lambda_i \right\} 
\end{equation}
for $u_1,\ldots  ,u_d\in \R^n\subset \t$ and $\lambda_1,\ldots ,\lambda_d\in \R$. 
Without loss of generality, we may assume that for $i=1,\ldots ,d$, each intersection 
$\Delta \cap \left\{  \xi \in \t^* \colon \< u_i, \xi \> =\lambda_i \right\}$ of $\Delta$ 
and the hyperplane defined by $\< u_i, \xi \> =\lambda_i$ is a {\it facet}, that is, a codimension one 
face of $\Delta$. We set 
\[
{\mathfrak J}=\left\{ I\subset \{1, \ldots, d\} \colon 
\bigcap_{i\in I}\{ \xi \in \Delta \colon \< u_i, \xi \> =\lambda_i \} \neq \emptyset 
\right\} .
\]
\begin{defn}\label{2.3}
The convex polytope $\Delta$ is said to be {\it Delzant}, if $\Delta$ satisfies the following 
conditions 
\begin{enumerate}
\item $\Delta$ is {\it rational}, that is, $u_1,\ldots  ,u_d\in \Z^n\subset \t$, 
\item $\Delta$ is {\it simple}, that is, exactly $n$ facets meet at each vertex of $\Delta$, 
\item $\Delta$ is {\it non-singular}, that is, $u_1,\ldots ,u_d$ (hence $\lambda_1, \ldots ,\lambda_d$) 
in (\ref{eq2.1}) can be taken so that $\left\{ u_i \right\}_{i\in I}$ is primitive for each 
$I\in {\mathfrak J}$. 
\end{enumerate}
\end{defn}
In the rest of this paper, when we say that a convex polytope $\Delta$ written as in 
$(\ref{eq2.1})$ is Delzant, we assume that $u_1,\ldots ,u_d$ are taken so that 
$\left\{ u_i \right\}_{i\in I}$ is primitive for each $I\in {\mathfrak J}$. 

Delzant shows that a symplectic toric manifold is exactly determined by the image of 
its moment map.  
\begin{thm}[\cite{D}]\label{2.4}
$(1)$ The image of a moment map of a symplectic toric manifold is a Delzant polytope. \\
$(2)$ By associating the image of a moment map to a symplectic toric manifold, 
the set of equivariantly symplectomorphism classes of $2n$-dimensional symplectic toric 
manifolds corresponds one-to-one to the set of Delzant polytopes in $\t^*\cong \R^n$ 
up to parallel transport in $\t^*$. 
\end{thm}

Theorem \ref{2.4} says that a symplectic toric manifold is recovered from a Delzant 
polytope. It is done as follows. Let $\Delta$ be an $n$-dimensional Delzant polytope 
defined by $(\ref{eq2.1})$. As in Example \ref{1.12}, for $i=1, \ldots ,d$, each vector 
$u_i$ in $(\ref{eq2.1})$ defines the circle actions on $(T^*T^n, \om_{T^*T^n})$ 
with the moment maps $\mu_{u_i} (\xi , \theta )=\< u_i, \xi \>$. 
Since these actions commute each other, we can obtain the simultaneous cut space 
$X_{\Delta}= \overline{(T^*T^n)}_{\{ \mu_{u_i}\ge \lambda_i \}_{i=1, \ldots ,d}}$ 
as we described in Remark \ref{1.9}. By the definition of the Delzant polytope and 
Theorem \ref{a3} in Appendix A, $X_{\Delta}$ is a $2n$-dimensional smooth symplectic 
manifold. Moreover by Remark \ref{1.10}, $X_{\Delta}$ is equipped with the Hamiltonian 
$T^n$-action with the moment map $\mu :X_{\Delta}\to \t^*$ which is induced from 
the natural $T^n$-action on $T^*T^n$ in Example \ref{1.3}. It is clear that this action is 
effective and the image of $\mu$ is $\Delta$. 
Hence the cut space $X_{\Delta}$ is the symplectic toric manifold which we want. 

From this construction, we can see that a symplectic toric manifold is locally identified with 
the Hamiltonian $T^n$-action on $\C^n$ in Example \ref{1.2} in the following sense. 
\begin{defn}
Let $\rho$ be an automorphism of $T^n$. Two $2n$-dimensional symplectic toric 
manifolds $(X_1, \om_1)$ and $(X_2, \om_2)$ are {\it $\rho$-equivariantly 
symplectomorphic}, if there exists a symplectomorphism 
$\varphi : (X_1, \om_1) \to (X_2, \om_2)$ such that $\varphi$ satisfies the two 
conditions 
\begin{enumerate}
\item $\varphi (t \cdot x)=\rho (t)\cdot \varphi (x)$ for $x\in X_1$ and 
$t \in T^n$, 
\item the following diagram is commutative     
\[
\xymatrix{
X_1\ar[d]_{\mu_1}\ar[r]^{\varphi}\ar@{}[dr]|\circlearrowright & X_2\ar[d]^{\mu_2} \\
\t^* & \t^*\ar[l]_{(d\rho )^*}, 
}
\]
where $\mu_1$ and $\mu_2$ are moment maps of $X_1$ and $X_2$, respectively.   
\end{enumerate}
\end{defn}
For a vertex $v\in \Delta$ which is defined by the exactly $n$ equalities 
$\displaystyle \< u_{i_a}, \xi \> =\lambda_{i_a}$ for $a=1, \ldots , n$, we define the 
open set $U_v\subset \Delta$ by 
\[
U_v=
\{ \xi \in \Delta \colon \< u_i, \xi \> > \lambda_i\ \text{for}\ i\neq i_1, 
\ldots , i_n   \} . 
\] 
Then $\{ \mu^{-1}(U_v) \}_{v:\text{vertex of }\Delta}$ is a open covering of 
$X_{\Delta}$. Moreover, for each vertex $v$, there exists an automorphism $\rho_v$ 
of $T^n$ such that $(d\rho_v^{-1})^*:\t^*\to \t^*$ sends $U_v$ diffeomorphically 
to the open set 
$(d\rho_v^{-1})^*(U_v)$ in $\t^*_{\ge 0}$ and $\mu^{-1}(U_v)$ is 
$\rho_v$-equivariantly symplectomorphic to  
$\mu^{-1}_{\C^n}((d\rho_v^{-1})^* (U_v))$ 
\[
\xymatrix{
X_{\Delta}\ar[d]_{\mu} \supset & \mu^{-1}(U_v)\ar[d]_{\mu} & \cong & 
\mu_{\C^n}^{-1}((d\rho_v^{-1})^*  (U_v))\ar[d]^{\mu_{\C^n}} & 
\subset \C^n\ar[d]^{\mu_{\C^n}} \\
\Delta \supset & U_v & \cong \ar@{}[u]|\circlearrowright  & 
(d\rho_v^{-1})^*(U_v) & {\subset} \t^*_{\ge 0}. \\
}
\] 
We shall show this claim. For a vertex $v$ which is defined as above, $\mu^{-1}(U_v)$ 
is identified with the open subset of the simultaneous cut space 
$\overline{(T^*T^n)}_{\{ \mu_{u_{i_a}}\ge \lambda_{i_a} \}_{a=1, \ldots ,n}}$ of 
the commutative Hamiltonian circle actions on $(T^*T^n)$ defined by $u_{i_j}$ with the 
moment maps $\mu_{u_{i_a}} (\xi , \theta )=\< u_{i_j}, \xi \>$ for $a=1, \ldots ,n$ 
naturally. Now define the parallel transport $p_v$ of $T^*T^n\times \C^k$ by 
\[
p_v(\xi , \theta , z) =(\xi -v, \theta , z) 
\] 
for $(\xi ,\theta , z)\in T^*T^n\times \C^k$. Then $p_v$ induces the equivariantly 
symplectomorphism 
from $\overline{(T^*T^n)}_{\{ \mu_{u_{i_a}}\ge \lambda_{i_a} \}_{a=1, \ldots ,n}}$ to 
$\overline{(T^*T^n)}_{\{ \mu_{u_{i_a}}\ge 0 \}_{a=1, \ldots ,n}}$. Then Theorem 
\ref{a4} in Appendix A with the above equivariantly symplectomorphisms implies the claim. 

\section{Twisted toric manifolds}
By the topological construction, a $2n$-dimensional symplectic toric manifold 
$X_{\Delta}$ is obtained from the trivial $T^n$-bundle on a Delzant polytope $\Delta$ 
by collapsing each fiber on the face $\{ \xi \in \Delta \colon \< u_i,\xi \> =\lambda_i \}$ 
by the circle subgroup $S^1_{u_i}$ generated by $u_i$. By replacing the trivial 
$T^n$-bundle on $\Delta$ to a $T^n$-bundle on an $n$-dimensional manifold 
with corners which may be non-trivial, we can obtain the notion of $2n$-dimensional 
twisted toric manifolds. 

\subsection{The definition and examples}
Let $B$ be an $n$-dimensional manifold with corners, $\pi_P :P\to B$ a principal 
$SL_n(\Z )$-bundle on $B$. The $T^n$-bundle and the $\Z^n$-bundle associated with $P$ 
by the natural action of $SL_n(\Z )$ on $T^n$ and on $\Z^n$ are denoted by 
$\pi_T:T^n_P\to B$ and $\pi_{\Z}:{\Z}^n_P\to B$, respectively. Consider a 
$2n$-dimensional manifold $X$, surjective maps $\nu :T^n_P\to X$ and $\mu :X\to B $ 
such that the following diagram is commutative 
\[
\xymatrix{
T^n_P \ar@{>>}[rr]^{\nu} \ar@{>>}[dr]_{\pi_T} &  & X \ar@{>>}[dl]^{\mu} \\
  & B. \ar@{}[u]|\circlearrowright &  
} 
\]
\begin{defn}\label{3.1}
The above tuple $\{ X, \nu ,\mu \}$ is called a $2n$-dimensional 
{\it twisted toric manifold} (or often called a {\it twisted toric structure} on $B$) 
associated with the principal 
$SL_n(\Z)$-bundle $\pi_P: P\to B$, if for arbitrary $b\in B$, there exist  
\begin{enumerate}
\item a coordinate neighborhood $(U, \varphi^B)$ of $b\in B$, that is, $U$ is an open 
neighborhood of $b$ in $B$ and $\varphi^B$ is an orientation consistent diffeomorphism 
from $U$ to the intersection $\t^*_{\ge 0}\cap D^n_{\epsilon}(\xi_0)$ of 
$\t^*_{\ge 0}$ 
and the open disc $D^n_{\epsilon}(\xi_0)$ in $\t^*=\R^n$ with a center 
$\xi_0\in \t^*_{\ge 0}$ and a radius $\epsilon >0$ which sends $b$ to $\xi_0$, 
(for the definition of $\t^*_{\ge 0}$, see Example \ref{1.2}, ) 
\item a local trivialization $\varphi^P:\pi_P^{-1}(U)\cong U\times SL_n(\Z )$ of $P$, 
(then $\varphi^P$ induces local trivializations 
$\varphi^T:\pi_T^{-1}(U)\cong U\times T^n$ and 
$\varphi^{\Z}:\pi_{\Z}^{-1}(U)\cong U\times \Z^n$ of $T^n_P$ and $\Z^n_P$, respectively, )
\item an orientation preserving diffeomorphism $\varphi^X:\mu^{-1}(U)\cong \mu_{\C^n}^{-1}(D^n_{\epsilon}(\xi_0))$ 
\end{enumerate}
such that the following diagram commutes
\[
\xymatrix{
\pi_T^{-1}(U)\ar[rrrr]^{\nu}\ar[drr]_{\pi_T}\ar[ddr]_{(\varphi^B\times \id_{T^n})\circ \varphi^T} & &  & & 
\mu^{-1}(U)\ar[dll]^{\mu}\ar[ddr]^{\varphi^X} & \\
  & & U\ar[ddr]^{\varphi^B}|\hole & & & \\
  & \varphi^B(U)\times T^n\ar[rrrr]^{\nu_{\C^n}}\ar[drr]_{\pr_1} & & & & 
\mu_{\C^n}^{-1}(D^n_{\epsilon}(\xi_0))\ar[dll]^{\mu_{\C^n}} \\
  &  & & \varphi^B(U), & &  
}
\]
where $\mu_{\C^n}$ is the moment map of $T^n$-action on $\C^n$ in Example \ref{1.2} 
and $\nu_{\C^n}$ is the map which is defined by 
\begin{equation}\label{eq3.1}
\nu_{\C^n}(\xi ,\theta )=(\sqrt{\xi_i}e^{2\pi \i \theta_i}). 
\end{equation}
Note that $\nu_{\C^n}$ is smooth only for $(\xi ,\theta )$ with all $\xi_i>0$. 
The tuple $(U, \varphi^P, \varphi^X, \varphi^B)$ is called a {\it locally toric chart}. 
If there are no confusions, we call simply $X$ a twisted toric manifold. 
\end{defn}
\begin{rem}[Orientations]\label{3.2}
We fix the orientations of $\t^*_{\ge 0}\times T^n$ (or $T^* T^n$) and $\C^n$ 
so that $d\xi_1\wedge d\theta_1 \wedge \cdots \wedge d\xi_n\wedge d\theta_n
(=(-1)^n\dfrac{(\om_{T^*T^n})^n}{n!})$ and $(-1)^n\dfrac{(\om_{\C^n})^n}{n!}$ 
are the positive volume forms, respectively. Then the map $\nu_{\C^n}$ in (\ref{eq3.1}) 
preserves the orientations. (c.f. Example \ref{1.11}) 
\end{rem}

\begin{ex}[Torus bundle]\label{ex3.3}
Let $\pi_P :P\to B$ be a principal $SL_n(\Z )$-bundle on a closed oriented 
$n$-dimensional manifold $B$. Then the associated $T^n$-bundle $\pi_T:T^n_P\to B$ 
itself is an example of a twisted toric manifold associated with $\pi_P:P\to B$. 
In particular, the even dimensional torus $T^{2n}$ is a twisted toric manifold, which is a  
$T^n$-bundle on $T^n$.  
\end{ex}
\begin{ex}[Symplectic toric manifold]\label{ex3.4}
A $2n$-dimensional symplectic toric manifold $X$ with a Delzant polytope $\Delta$ has 
a structure of a twisted toric manifold associated with the trivial $SL_n(\Z )$-bundle 
on $\Delta$. In fact, as we described in Section 3, $X$ is obtained from the trivial 
$T^n$-bundle on $\Delta$ by collapsing each fiber on the facet of $\Delta$ by the circle 
subgroup which is generated by the inward pointing normal vector in $\Z^n$ of the facet.  
\end{ex}
\begin{ex}[$S^{2n-1}\times T^{2n-3}$]\label{ex3.5}
For $n\ge 2$, let $S^{2n-1}$ be the unit sphere in $\C^n$. The product 
$X=S^{2n-1}\times T^{2n-3}$ is a twisted toric manifold associated with 
the trivial $SL_{2n-2}(\Z )$-bundle $P$ on the $(2n-2)$-dimensional unit disk 
$B =\overline{D^{2n-2}}=
\{ z\in \C^{n-1} \colon \lVert z\rVert \le 1 \}$ in $\C^{n-1}$. The maps 
$\nu :T^{2n-2}_P=B\times T^{2n-2}\to X$ and 
$\mu : X\to B$ are given by 
\[
\begin{split}
\nu (z, \theta )&=\left(
(\sqrt{1-\lVert z\rVert^2}e^{2\pi \i \theta_1}, z), (\theta_2, \ldots , \theta_{2n-2} )
\right) ,\\
\mu (w, \tau)&=(w_2,\ldots , w_{n-1})
\end{split}
\] 
for $(z, \theta )\in B\times T^{2n-2}$, and $(w, \tau )\in S^{2n-1}\times T^{2n-3}$.  
\end{ex}
\begin{ex}[Even dimensional sphere]\label{ex3.6}
Let $\Delta^n$ be the $n$-simplex and $\partial \Delta^n_0$ the facet of $\Delta^n$ which are 
defined by 
\[
\begin{split}
\Delta^n &=\{ \xi =(\xi_i) \in \t^*=\R^n \colon \xi_i\ge 0,\ \sum_{i=1}^n\xi_i\le 1\} ,\\
\partial \Delta^n_0&=\{ \xi =(\xi_i) \in \Delta^n \colon \sum_{i=1}^n\xi_i=1\} . 
\end{split}
\]
The $2n$-dimensional sphere $S^{2n}$ in $\C^n\times \R$ is equipped with a twisted toric structure 
associated with the trivial $SL_n(\Z )$-bundle $P$ on the quotient space 
$B=\Delta^n/\partial \Delta^n_0$ of $\Delta^n$ by $\partial (\Delta^n)_0$ which is explained 
as follows. $S^{2n}$ can be thought of as the one point compactification 
$\overline{D^{2n}}/\{ z\in \overline{D^{2n}} \colon \lVert z\rVert =1\}$ of the $2n$-dimensional 
unit disk $\overline{D^{2n}}=\{ z\in \C^n \colon \lVert z\rVert \le 1\}$ in $\C^n$. 
With this identification, the maps $\overline{\nu}:\Delta^n\times T^n\to \overline{D^{2n}}$ and 
$\overline{\mu}:\overline{D^{2n}}\to \Delta^n$ defined by
\[
\begin{split}
\overline{\nu}(\xi , \theta )&=(\sqrt{\xi_i}e^{2\pi \i \theta_i}) ,\\
\overline{\mu}(z)&=(\lvert z_i\rvert^2)
\end{split}
\]
for $(\xi , \theta )\in \Delta^n\times T^n$ and $z=(z_i)\in \overline{D^{2n}}$ induce the maps 
$\nu :T^n_P=B \times T^n\to S^{2n}$ and $\mu :S^{2n}\to B$, respectively.  
\end{ex}
\begin{ex}\label{ex3.7} 
Consider the Hamiltonian $T^2$-action on 
$(T^*T^2\times \C^2, \om_{T^*T^2}\oplus \om_{\C^2})$ defined by 
\[
t\cdot (\xi ,\theta ,z)=
(\xi , (\theta_1, \theta_2+t_1-t_2), (e^{-2\pi \i t_1}z_1, e^{-2\pi \i t_2}z_2)) 
\]
with the moment map 
$\Phi :T^*T^2\times \C^2\to \R^2$
\[
\Phi (\xi , \theta ,z)=
(\xi_2-\lvert z_1\rvert^2, -\xi_2-\lvert z_2\rvert^2+1 ). 
\]
We should remark that the second component of $\Phi$ is added the constant 
$1$. (cf. the end of Definition \ref{1.1}.) We denote by $\overline{X}$ its symplectic 
quotient $\Phi^{-1}(0)/T^2$.  
Define the right action of $\Z$ on $T^*T^2\times \C^2$ by
\[
(\xi , \theta ,z)\cdot n =((\xi_1+n, \xi_2), \rho (-n)\theta , \varphi_n(z))
\]
for $(\xi , \theta ,z)\in T^*T^2\times \C^2$ and $n\in \Z$,  
where $\rho :\Z\to SL_2(\Z )$ is the homomorphism 
\[
\rho (n)=
\begin{pmatrix}
-1 & 0 \\
0 & -1
\end{pmatrix}^n
\]
and 
\[
\varphi_n (z)=
\begin{cases}
z & n:\text{even} \\
\bar{z} & n:\text{odd}
\end{cases} . 
\]
It is easy to see that the action descends to the action of $\Z$ on 
$\overline{X}$, and its quotient space $\overline{X}/\Z$ is denoted by $X$. 
We shall explain that $X$ is a twisted toric manifold associated with 
the principal $SL_2(\Z )$-bundle $P$ on the cylinder $B=S^1\times [0,1]$ 
which is determined by the representation 
$\rho :\Z =\pi_1(B)\to SL_2(\Z )$. 
Let $\widetilde{B}=\R \times [0,1]$ be the universal covering of 
$B$ on which the fundamental group $\pi_1(B)=\Z$ acts as a 
deck transformation by 
\[
\xi \cdot n=(\xi_1+n, \xi_2). 
\] 
Since $\overline{X}$ is a simultaneous cut space by two $S^1$-actions on 
$T^*T^2$ corresponding to the positive and negative second fundamental 
vectors $e_2=(0,1)$, $-e_2=(0,-1)$ in $\Z^2$, the natural $T^2$-action on $T^*T^2$ 
induces the Hamiltonian action of $T^2$ on $\overline{X}$ with the 
moment map $\overline{\mu}([\xi , \theta ,z])=\xi $ whose image is equal to 
$\widetilde{B}$. Consider the following commutative triangle of maps 
\begin{equation}\label{triangle}
\xymatrix{
\widetilde{B}\times T^2\ar[rr]^{\overline{\nu}} \ar[dr]_{\pr_1} & 
\ar@{}[d]|\circlearrowright & \overline{X}\ar[dl]^{\overline{\mu}} \\
  & \widetilde{B}, & 
}
\end{equation}
where $\overline{\nu}:\widetilde{B}\times T^2\to \overline{X}$ is given by
\[
\overline{\nu}(\xi , \theta )=[\xi ,\theta , (\sqrt{\xi_2}, \sqrt{1-\xi_2})]. 
\]
$\Z$ acts on $\widetilde{B}\times T^2$ by 
\[
(\xi , \theta )\cdot n=((\xi_1+n, \xi_2), \rho(-n)\theta ) 
\]
for $(\xi , \theta )\in \widetilde{B}\times T^2$ and $n\in \Z$. 
Since the maps $\pr_1$, $\overline{\nu}$, and $\overline{\mu}$ are equivariant with 
respect to the above $\Z$-actions, these descend to the maps $\pi_T:T^2_P\to B$, 
$\nu :T^2_P\to X$, and $\mu :X\to B$. 
\end{ex}
\begin{ex}\label{ex3.8}
Let $B$ be a compact, connected, and oriented surface of genus $g\ge 1$ 
with one boundary component and $k(\neq 1)$ corner points, $B_1$ the open 
subset of $B$ which is obtained by removing a sufficiently small closed 
neighborhood of the boundary from $B$, and $B_2$ an open neighborhood of 
the boundary such that $B_1 \cap B_2 \cong S^1\times [0,1]$, see Figure \ref{fig2}. $B$ has the oriented boundary loop which we denote by $\gamma$. 
Moreover, for $k>0$, the boundary $\partial B$ of $B$ consists of exactly 
$k$ edge arcs which are denoted by $\gamma_1$, $\ldots$, $\gamma_k$ as in 
Figure \ref{fig2}.
\begin{figure}[htbp]
\begin{center}
\input{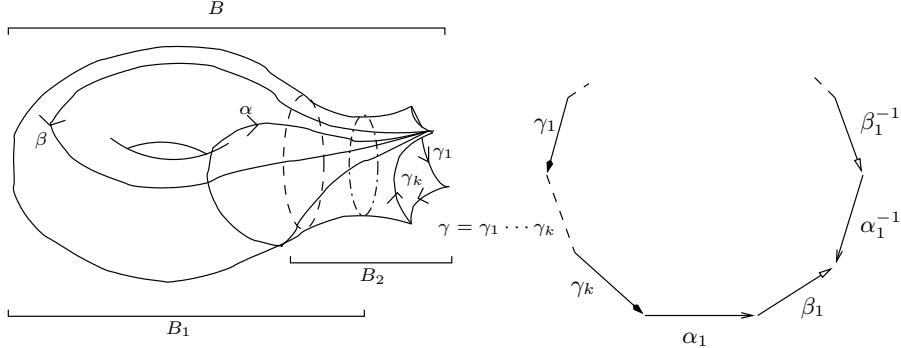}
\caption{$B$, $B_1$, and $B_2$}
\label{fig2}
\end{center}
\end{figure}
The fundamental group $\pi_1(B )$ of $B$ can be identified with the quotient 
group $F_{2g+1}/N$ of the free group $F_{2g+1}$ of $2g+1$ generators $\alpha_1$, 
$\beta_1$, $\ldots$, $\alpha_g$, $\beta_g$, $\gamma$ by the least normal subgroup 
$N$ containing $\prod_{i=1}^g[\alpha_i, \beta_i]\gamma$, and define the representation 
$\rho :\pi_1(B )\to SL_2(\Z )$ by 
\[
\rho (\alpha_i )=
\begin{pmatrix}
1 & a_i \\
0 & 1
\end{pmatrix},\ 
\rho (\beta_i )=
\begin{pmatrix}
1 & b_i \\
0 & 1
\end{pmatrix},\ 
\rho (\gamma )=
\begin{pmatrix}
1 & 0 \\
0 & 1
\end{pmatrix}.  
\]
$\rho$ determines the principal $SL_2(\Z )$-bundle $P$ on $B$ as 
the quotient space $\widetilde{B}\times_{\rho}SL_n(\Z )$ of the 
action of $\pi_1(B )$ on the product $\widetilde{B}\times SL_n(\Z )$ 
of the universal cover $\widetilde{B}$ of $B$ and $SL_n(\Z )$ 
which is defined by 
\[
(\widetilde{b}, g)\cdot \tau =(\widetilde{b}\cdot \tau , \rho(\tau^{-1})g)
\] 
for $(\widetilde{b}, g)\in \widetilde{B}\times SL_n(\Z )$ and $\tau \in \pi_1(B)$, 
where $\widetilde{b}\cdot \tau$ means the deck transformation. 
Its associated $T^2$-bundle by the natural action of 
$SL_2(\Z )$ on $T^2$ is denoted by $T^2_P$. 

For $k\ge 3$, we fix the 4-dimensional symplectic 
toric manifold $X_{\Delta}$ corresponding to the Delzant polytope $\Delta$ with 
$k$ vertices, and let $X_2$ denote the subspace of $X_{\Delta}$ which is obtained by 
removing from $X_{\Delta}$ the inverse image $\mu^{-1}(\overline{D^2})$ of a 
small closed disk $\overline{D^2}$ 
in the interior of $\Delta$ by the moment map $\mu :X_{\Delta}\to \t^*=\R^2$. 
We take a small open disk ${D'}^2$ including $\overline{D^2}$ in the interior of $\Delta$. 
Since $T^2_P\big|_{B_1 \cap B_2}\to B_1 \cap B_2$ is trivial, it can be identified with 
$\mu^{-1}({D'}^2\backslash \overline{D^2}) \to {D'}^2\backslash \overline{D^2}$. 
Then we can glue $T^2_P\big|_{B_1}\to B_1$ with 
$X_{\Delta}\backslash \mu^{-1}(\overline{D^2})\to \Delta \backslash \overline{D^2}$ by 
this identification to obtain a $4$-dimensional twisted toric manifold associated 
with $P$.  

For $k=0$ (resp. $k=2$), by replacing $X_{\Delta}$ by $S^3\times S^1$ 
in Example \ref{ex3.5} (resp. $S^4$ in Example \ref{ex3.6}) in the above 
construction, we can obtain a $4$-dimensional twisted toric manifold 
in this case. 
\begin{figure}[hbtp]
\begin{center}
\input{fig3.pstex_t}
\caption{gluing $B_1$ and $\Delta \backslash \overline{D^2}$}
\label{fig3}
\end{center}
\end{figure}
\end{ex}
\begin{ex}\label{ex3.9}
Let $B$ be a compact, connected, and oriented surface of genus one 
with one boundary component and one corner points, $B_1$ the set 
of interior points of $B$, that is, $B_1=B\backslash \partial B$ . 
In this case, consider the principal $SL_2(\Z )$-bundle $P$ on $B$ which is 
determined by the representation $\rho :\pi_1(B )\to SL_2(\Z )$
\[
\rho (\alpha )=
\begin{pmatrix}
1 & 0 \\
-1 & 1
\end{pmatrix},\ 
\rho (\beta )=
\begin{pmatrix}
1 & -1 \\
0 & 1
\end{pmatrix},\ 
\rho (\gamma )=
\begin{pmatrix}
3 & 1 \\
-1 & 0
\end{pmatrix}.  
\]
As above, we also denote by $\pi_T:T^2_P\to B$ its associated 
$T^2$-bundle by the natural action of $SL_2(\Z )$ on $T^2$. 
In this case, $\pi_T:T^2_P\to B$ is no more trivial near the boundary 
$\partial B$ of $B$, and let us construct the twisted toric structure near 
$\partial B$ as follows. We define the subset $\overline{B}_2$ of $\t^*_{\ge 0}$ 
by 
\[ 
\overline{B}_2=\{ \xi \in \t^*(=\R^2) \colon 0\le \xi_1 <4,\ 0\le \xi_2<1\} 
\cup \{ \xi \in \t^* \colon 0\le \xi_1 <1,\ 0\le \xi_2<4\} . 
\]
The restrictions $\nu_{\C^2}\big|_{\overline{B}_2\times T^2}:
\overline{B}_2\times T^2\to \mu_{\C^2}^{-1}(\overline{B}_2)$ 
and $\mu_{\C^2}\big|_{\mu_{\C^2}^{-1}(\overline{B}_2)}:
\mu_{\C^2}^{-1}(\overline{B}_2)\to \overline{B}_2$ of 
$\nu_{\C^2}$ and $\mu_{\C^2}$ defined in (\ref{eq3.1}) and 
Example \ref{1.2}, respectively, form a commutative triangle of surjective maps 
\begin{equation}\label{tri}
\xymatrix{
\overline{B}_2\times T^2 \ar@{>>}[rr]^{\nu_{\C^2}} \ar@{>>}[dr]_{\pr_1} &  & 
\mu_{\C^2}^{-1}(\overline{B}_2) \ar@{>>}[dl]^{\mu_{\C^2}} \\
  & \overline{B}_2. \ar@{}[u]|\circlearrowright &  
} 
\end{equation}
Let $U_1$ and $U_2$ be open sets of $\overline{B}_2$ which are defined by
\[
\begin{split}
U_1&=\{ \xi \in \t^* \colon 3<\xi_1<4,\ 0\le \xi_2<1\} ,\\
U_2&=\{ \xi \in \t^* \colon 0\le \xi_1<1,\ 3<\xi_2<4\} . 
\end{split} 
\]
We define diffeomorphisms $\varphi^B :U_1\to U_2$, 
$\varphi^T:U_1\times T^2\to U_2\times T^2$, and 
$\varphi^X :\mu_{\C^2}^{-1}(U_1)\to \mu_{\C^2}^{-1}(U_2)$ by 
\[
\begin{split}
& \varphi^B (\xi )=(\xi_2, 7-\xi_1), \\
& \varphi^T (\xi ,\theta )=(\varphi^B(\xi ), \rho (\gamma )\theta ), \\
& \varphi^X (z)=(\dfrac{z_1^3z_2}{\abs{z_1}^3}, 
             \sqrt{7-\abs{z_1}^2}\left( \dfrac{z_1}{\abs{z_1}}\right)^{-1}). 
\end{split}
\]
We denote by $X_2$ the manifold which is obtained from 
$\mu_{\C^2}^{-1}(\overline{B}_2)$ by gluing $\mu_{\C^2}^{-1}(U_1)$ and 
$\mu_{\C^2}^{-1}(U_2)$ with $\varphi^X$ and denote by $B_2$ the surface 
with one corner which is obtained from $\overline{B}_2$ by gluing $U_1$ and 
$U_2$ with $\varphi^B$. 
\begin{figure}[hbtp]
\begin{center}
\input{fig4.pstex_t}
\caption{$B_2$}
\label{fig4}
\end{center}
\end{figure}
$B_2$ can be naturally identified with a neighborhood of the boundary of $B$. 
Since the following diagram commutes 
\[
\xymatrix{
U_1\times T^2\ar[rrr]^{\varphi^T}\ar[drr]^{\nu_{\C^2}}\ar[ddr]_{\pr_1} &  
&  & U_2\times T^2\ar[drr]^{\nu_{\C^2}}\ar[ddr]|\hole ^{\pr_1} 
&  &  \\
  &  & \mu_{\C^2}^{-1}(U_1)\ar[rrr]^{\varphi^X\hspace*{1.4cm}}\ar[dl]^{\mu_{\C^2}} &  
  &  & \mu_{\C^2}^{-1}(U_2)\ar[dl]^{\mu_{\C^2}} \\
  & U_1 \ar[rrr]^{\varphi^B} &  &  & U_2, &  
} 
\]
under the identification of $B_2$ with a neighborhood of $\partial B$, 
the restriction $\pi_T\big|_{B_2}:T^2_P\big|_{B_2}\to B_2$ of the torus bundle 
$\pi_T:T^2_P\to B$ to $B_2$ is obtained from the trivial bundle 
$\pr_1:\overline{B}_2\times T^2\to \overline{B}_2$ by gluing  $U_1\times T^2$ and 
$U_2\times T^2$ by $\varphi^T$, and the diagram (\ref{tri}) induces the following 
commutative triangle of surjective maps 
\[
\xymatrix{
T^2_P\big|_{B_2} \ar@{>>}[rr]^{\nu_2} \ar@{>>}[dr]_{\pi_T} &  & 
X_2 \ar@{>>}[dl]^{\mu_2} \\
  & B_2, \ar@{}[u]|\circlearrowright &  
} 
\]
where $\nu_2$ and $\mu_2$ are maps induced by $\nu_{\C^2}$ and $\mu_{\C^2}$, 
respectively. It is easy to see that the restriction $\mu_2 \big|_{\mu_2^{-1}(B_1\cap B_2)}:
\mu_2^{-1}(B_1\cap B_2)\to B_1\cap B_2$ is a $T^2$-bundle with the structure group 
$SL_2(\Z )$ and the restriction of $\nu_2$ to $T^2_P\big|_{B_1\cap B_2}$ 
is a bundle isomorphism from $T^2_P\big|_{B_1\cap B_2}$ to $\mu_2^{-1}(B_1\cap B_2)$. 
Thus we can patch $\pi_T\big|_{B_1}:T^2_P\big|_{B_1}\to B_1$ with  $\mu_2 :X_2\to B_2$ 
by this isomorphism to get the twisted toric manifold $\mu :X\to B$ associated 
with $\pi :P\to B$. 
\end{ex}

\subsection{A remark on compatible symplectic forms} 
Let $B$ be an $n$-dimensional manifold and $X$ a $2n$-dimensional twisted toric manifold 
associated with a principal $SL_n(\Z )$-bundle $P\to B$ on $B$. For each locally toric chart 
$( U, \varphi^P, \varphi^X, \varphi^{B})$ of $X$, we have a symplectic form on 
$\mu^{-1}(U)$ which comes from $\om_{\C^n}$ on $(\mu_{\C^n}^{-1})(\varphi^{B}(U))$. 
In the definition of twisted toric manifolds, we do not assume that these local symplectic forms 
are patched together to a global symplectic form of $X$.  
In this section, we shall investigate the condition for these forms to be patched together 
to a global symplectic form on $X$. 
Let $(U_{\alpha}, \varphi^P_{\alpha}, \varphi^X_{\alpha}, \varphi^{B}_{\alpha})$, 
$(U_{\beta}, \varphi^P_{\beta}, \varphi^X_{\beta}, \varphi^{B}_{\beta})$ 
be two locally toric charts of $X$ which satisfies 
$U_{\alpha \beta}=U_{\alpha}\cap U_{\beta}\neq \emptyset$. 
We denote by $g_{\beta \alpha}: U_{\alpha \beta}\to SL_n(\Z)$ the transition 
function of $P$ and set 
\[
\begin{split}
\varphi^{B}_{\beta \alpha}&=
(\varphi^{B}_{\beta}\big|_{U_{\alpha \beta}})\circ 
(\varphi^{B}_{\alpha}\big|_{U_{\alpha \beta}})^{-1}:
\varphi^{B}_{\alpha}(U_{\alpha \beta})\to \varphi^{B}_{\beta}(U_{\alpha \beta}), \\
\varphi^T_{\beta \alpha}&=\varphi^T_{\beta}\big|_{\pi_T^{-1}(U_{\beta \alpha})}
\circ (\varphi^T_{\alpha}\big|_{\pi_T^{-1}(U_{\beta \alpha})})^{-1}: 
U_{\alpha \beta}\times T^n \to U_{\alpha \beta}\times T^n, \\
\varphi^X_{\beta \alpha}&=(\varphi^X_{\beta}\big|_{\mu^{-1}(U_{\alpha \beta})})\circ 
(\varphi^X_{\alpha}\big|_{\mu^{-1}(U_{\alpha \beta})})^{-1}:
(\mu_{\C^n}^{-1}(\varphi^{B}_{\alpha}(U_{\alpha \beta}))) \to 
(\mu_{\C^n}^{-1}(\varphi^{B}_{\beta}(U_{\alpha \beta}))). 
\end{split}
\] 
We also define the  maps 
$\widetilde{g}_{\beta \alpha}:\varphi^{B}_{\alpha}(U_{\alpha \beta})\to SL_n(\Z)$ 
and 
$\widetilde{\varphi}^T_{\beta \alpha}:\varphi^{B}_{\alpha}(U_{\alpha \beta})\times T^n
\cong \varphi^{B}_{\beta}(U_{\alpha \beta})\times T^n$ by 
\[
\begin{split}
\widetilde{g}_{\beta \alpha}&=
g_{\beta \alpha}\circ (\varphi^{B}_{\alpha}\big|_{U_{\alpha \beta}})^{-1}, \\
\widetilde{\varphi}^T_{\beta \alpha}&=(\varphi^{B}_{\beta }\big|_{U_{\alpha \beta}}
\times \id_{T^n})\circ \varphi^T_{\beta \alpha}\circ 
(\varphi^{B}_{\alpha}\big|_{U_{\alpha \beta}}\times \id_{T^n})^{-1}. 
\end{split}
\]
Note that they satisfy the following conditions 
\[
\begin{split}
\varphi^T_{\beta \alpha}(b, \theta )&=(b, g_{\beta \alpha}(b)\theta ), \\
\widetilde{\varphi}^T_{\beta \alpha}(\xi , \theta' )&=
(\varphi^B_{\beta \alpha}(\xi ),\widetilde{g}_{\beta \alpha}(\xi )\theta' )
\end{split}
\]
for $(b,\theta )\in U_{\beta \alpha}\times T^n$ and 
$(\xi ,\theta' ) \in \varphi^B_{\alpha}(U_{\alpha \beta})\times T^n$ 
\[
\xymatrix{
   & \pi_T^{-1}(U_{\alpha \beta})\ar[dl]_{(\varphi^{B}_{\alpha}\times \id_{T^n})\circ \varphi^T_{\alpha}} \ar[dr]^{(\varphi^{B}_{\beta}\times \id_{T^n})\circ \varphi^T_{\beta}} 
\ar@{-}[d]^{\nu} & \\
\varphi^B_{\alpha}(U_{\alpha \beta})\times T^n\ar[rr]^{\widetilde{\varphi}^T_{\beta \alpha}\ \ \ \ } \ar[dd]^{\nu_{\C^n}} & \ar[d] & \varphi^B_{\beta}(U_{\alpha \beta})\times T^n
\ar[dd]^{\nu_{\C^n}} \\
& \mu^{-1}(U_{\alpha \beta})\ar[dl]_{\varphi^X_{\alpha}} \ar[dr]^{\varphi^X_{\beta}} 
\ar@{-}[d]^{\mu} &  \\
\mu^{-1}_{\C^n}(\varphi^B_{\alpha}(U_{\alpha \beta}))\ar[rr]^{\varphi^X_{\beta \alpha}\ \ \ \ } \ar[dd]^{\mu_{\C^n}} & \ar[d] & \mu^{-1}_{\C^n}(\varphi^B_{\beta}(U_{\alpha \beta}))
\ar[dd]^{\mu_{\C^n}} \\
 &  U_{\alpha \beta} \ar[dl]_{\varphi^B_{\alpha}} \ar[dr]^{\varphi^B_{\beta}}& \\ 
\varphi^B_{\alpha}(U_{\alpha \beta})\ar[rr]^{\varphi^B_{\beta \alpha}\ \ \ \ } & & \varphi^B_{\beta}(U_{\alpha \beta}).
}
\]
\begin{prop}
$\varphi^X_{\beta \alpha}$ preserves the symplectic form $\om_{\C^n}$ on 
$\mu^{-1}_{\C^n}(\varphi^B_{\alpha}(U_{\alpha \beta}))$, if and only if, 
up to additive constant, $\varphi^{B}_{\beta \alpha}$ is of the form 
$\varphi^{B}_{\beta \alpha}(\xi )= ^t(\widetilde{g}_{\beta \alpha}(\xi ))^{-1}\xi$. 
\end{prop}
\begin{proof}
Since for $k=\alpha$, $\beta$, the symplectic form $\om_{\C^n}$ on 
$(\mu_{\C^n}^{-1}(\varphi^B_k(U_k)))$ comes from $\om_{T^*T^n}$ on 
$\varphi^B_k(U_k)\times T^n\subset T^*T^n$ by symplectic cutting, it is 
sufficient to see the condition that $\widetilde{\varphi}^T_{\beta \alpha}$ 
preserves $\om_{T^*T^n}$. Whereas 
\[
( \widetilde{\varphi}^T_{\beta \alpha})^* \om_{T^*T^n}=\sum_{j,k=1}^n
\left( \sum_{i=1}^n(\widetilde{g}_{\beta \alpha}(\xi ))_{ij}
\dfrac{\partial (\varphi^{B}_{\beta \alpha})_i}{\partial \xi_k}\right)
d\theta_j\wedge dx_k,  
\] 
$\widetilde{\varphi}^T_{\beta \alpha}$ preserves $\om_{T^*T^n}$, if and only if 
\[
\sum_{i=1}^n(\widetilde{g}_{\beta \alpha}(\xi ))_{ij}
\dfrac{\partial (\varphi^{B}_{\beta \alpha})_i}{\partial \xi_k}=\delta_{jk},  
\]
that is, the Jacobi matrix 
$(\dfrac{\partial (\varphi^{B}_{\beta \alpha})_i}{\partial \xi_j})$ is equal to 
$^t(\widetilde{g}_{\beta \alpha}(\xi ))^{-1}$. Since $\widetilde{g}_{\beta \alpha}$ 
is locally constant, this implies the lemma.  
\end{proof}
\begin{defn}
Let $\om$ be a symplectic form on a twisted toric manifold $X$. 
$\om$ is said to be {\it compatible} with respect to the twisted toric 
structure of $X$, if for each locally toric chart $( U, \varphi^P, \varphi^X, \varphi^{B})$, 
the restriction of $\om$ to $\mu^{-1}(U)$ is equal to ${\varphi^X}^*\om_{\C^n}$. 
\end{defn}
\begin{ex}
A $2n$-dimensional torus $T^{2n}$ is a twisted toric manifold, because $T^{2n}$ itself is 
the trivial $n$-dimensional torus bundle on a $n$-dimensional torus. The symplectic 
structure on $T^{2n}$ which is induced from that of $\R^{2n}$ is compatible.  
\end{ex}
\begin{ex}
The symplectic structure of every symplectic toric manifold is compatible with  respect to 
the natural twisted toric structure.   
\end{ex}
\begin{ex}\label{ex5.5}
The product of an even dimensional torus and a symplectic toric manifold has a compatible 
symplectic structure. 
\end{ex}
\begin{ques}
Are there any twisted toric manifold associated with a non trivial principal $SL_n(\Z)$-bundle 
which has a compatible symplectic form$?$
\end{ques}

\section{The classification theorem}
In this section, we shall prove the classification theorem for twisted toric manifolds. 
Let $B$ be an $n$-dimensional manifold and $\pi_P:P\to B$ a principal 
$SL_n(\Z )$-bundle $P\to B$ on $B$. We continue to use same notations for 
associated $T^n$- and $\Z^n$-bundles as in Section 4. In this section, we assume 
that $\partial B\neq \emptyset$. For arbitrary $b\in B$, define 
\[
n (b)=\# \{ i \colon \varphi^{B}(b)_i=\< e_i, \varphi^{B}(b) \> =0\} ,
\]
where $e_i$ is the $i$th fundamental vector 
$e_i=^t\stackrel{i}{(0,\ldots ,0,1,0,\ldots ,0)}$ of $\Z^n$ and 
$\varphi^{B}$ is a local coordinate on a neighborhood of $b$ defined as in (i) of Definition 
\ref{3.1}. Note that $n(b)$ does not depend on the choice of $\varphi^{B}$. 
Let ${\mathcal S}^{(k)}B$ be the $k$-dimensional strata of $B$ which is defined by
\[
{\mathcal S}^{(k)}B =\{ b\in B \colon n(b)=n-k \} . 
\]
Let $\pi_{\mathcal L}:{\mathcal L}\to {\mathcal S}^{(n-1)}B$ be a rank one sub-lattice 
bundle of the restriction 
$\pi_{\Z}\big|_{{\mathcal S}^{(n-1)}B}:\Z^n_P\big|_{{\mathcal S}^{(n-1)}B}
\to {\mathcal S}^{(n-1)}B$ of the associated $\Z^n$-bundle $\pi_{\Z}:\Z^n_P\to B$ 
to the codimension one strata ${\mathcal S}^{(n-1)}B$. 
\begin{defn}\label{4.2}
$\pi_{\CL}:\CL \to \CS^{(n-1)}B$ is said to be {\it primitive}, if for arbitrary point $b$ 
of the $k$-dimensional strata ${\mathcal S}^{(k)}B$, there exist
\begin{enumerate}
\item an open neighborhood $U$ of $b\in B$ 
whose intersection $U\cap \CS^{(n-1)}B$ 
with $\CS^{(n-1)}B$ has exactly $n-k$ connected components $(U\cap {\mathcal S}^{(n-1)}B)_1$, 
$\ldots$, $(U\cap {\mathcal S}^{(n-1)}B)_{n-k}$, 
\item a local trivialization $\varphi^P:P\big|_U\cong U\times SL_n(\Z )$, 
\item a primitive tuple $\{ L_1, \ldots ,L_{n-k} \}$ of rank one sub-lattices in $\Z^n$ 
\end{enumerate}
such that the associated local trivialization $\varphi^{\Z}:\Z^n_P\big|_U\cong U\times \Z^n$ of 
$\varphi^P$ maps each ${\mathcal L}\big|_{(U\cap {\mathcal S}^{(n-1)}B)_a}$ to 
$(U\cap {\mathcal S}^{(n-1)}B)_a\times L_a$ 
\[
\begin{matrix}
\Z^n_P\big|_U & \overset{\varphi^{\Z}}{\cong} & U\times \Z^n \\
\cup          &       & \cup \\
\Z^n_P\big|_{(U\cap {\mathcal S}^{(n-1)}B)_a} & \cong & 
(U\cap {\mathcal S}^{(n-1)}B)_a\times \Z^n \\
\cup          &       & \cup \\
{\mathcal L}\big|_{(U\cap {\mathcal S}^{(n-1)}B)_a} & \cong & 
(U\cap {\mathcal S}^{(n-1)}B)_a\times L_a
\end{matrix}
\]
for $a=1, \ldots ,n-k$. 
\end{defn}
\begin{rem}
(1) Definition \ref{4.2} does not depend on the choice of a neighborhood $U$ in (i) of 
Definition \ref{4.2}, since 
the notion of primitivity is invariant under the action of $SL_n(\Z )$. (cf. Remark \ref{2.25})\\
(2) The automorphism group $\Aut (P)$ of $P$ acts on the set of primitive rank one 
sub-lattice bundles of $\pi_{\Z}:\Z^n_P\bigl|_{\CS^{(n-1)}B}\to \CS^{(n-1)}B$ as the 
automorphisms of the restriction of the associated lattice bundle $\pi_{\Z}:\Z^n_P\to B$ 
to $\CS^{(n-1)}B$.
\end{rem}
\begin{thm}\label{thm5.3}
For any twisted toric manifold $X$ associated with a principal $SL_n(\Z )$-bundle $\pi_P:P\to B$, 
there exists a primitive rank one sub-lattice bundle $\pi_{\CL_X}:\CL_X \to \CS^{(n-1)}B$ of 
$\pi_{\Z}:\Z^n_P\bigl|_{\CS^{(n-1)}B}\to \CS^{(n-1)}B$ which is determined uniquely by 
$X$. $\pi_{\CL_X}:\CL_X \to \CS^{(n-1)}B$ is called a {\it characteristic bundle} of $X$. 
\end{thm}
To prove Theorem \ref{thm5.3}, we need some preliminaries. 
Let $(U_{\alpha}, \varphi^P_{\alpha}, \varphi^X_{\alpha}, \varphi^{B}_{\alpha})$, 
$(U_{\beta}, \varphi^P_{\beta}, \varphi^X_{\beta}, \varphi^{B}_{\beta})$ be two 
locally toric charts of $X$ with non empty intersection 
$U_{\alpha \beta}=U_{\alpha}\cap U_{\beta}\neq \emptyset$. We also use the notations 
defined in Section 4.2. Let $\Z_i$ be the rank one sub-lattice of the free $\Z$-module $\Z^n$ 
which is spanned by $e_i$. For $j=\alpha$, $\beta$, we define 
\[
{\mathfrak I}_j=
\left\{  i\in\{ 1, \ldots ,n \} \colon \{ \xi \in \varphi^B_j(U_{\alpha \beta}) \colon \< \xi ,e_i\> =\xi_i=0 \} \neq \emptyset \right\} . 
\]
For any $i_{\alpha}\in {\mathfrak I}_{\alpha}$, there exists a unique 
$i_{\beta}\in {\mathfrak I}_{\beta}$ such that 
\[
\varphi^{B}_{\beta \alpha}(\{ \xi \in \varphi^B_{\alpha}(U_{\alpha \beta}) \colon \< \xi ,e_{i_{\alpha}}\> =\xi_{i_{\alpha}}=0 \} ) 
= \{ \xi \in \varphi^B_{\beta}(U_{\alpha \beta}) \colon \< \xi ,e_{i_{\beta}}\> =\xi_{i_{\beta}}=0 \} . 
\] 
\begin{lem}\label{4.1}
For such $i_{\alpha}\in {\mathfrak I}_{\alpha}$ and $i_{\beta}\in {\mathfrak I}_{\beta}$, 
\[
g_{\beta \alpha}(b)\Z_{i_{\alpha}}=\Z_{i_{\beta}} 
\] 
for all $b\in (\varphi^{B}_{\alpha})^{-1}
(\{ \xi \in \varphi^B_{\alpha}(U_{\alpha \beta}) \colon \< \xi ,e_{i_{\alpha}}\> =\xi_{i_{\alpha}}=0 \} )
\subset U_{\alpha \beta}$. 
\end{lem}
\begin{proof}
Since $g_{\beta \alpha}$ is locally constant, it is sufficient to show 
the lemma on the set 
\[
(\varphi^{B}_{\alpha})^{-1}(\{ \xi \in \varphi^B_{\alpha}(U_{\alpha \beta}) \colon 
\xi_{i_{\alpha}}=0,\ \xi_i>0,\ \text{for}\ i\neq i_{\alpha} \} ). 
\]
Let $b\in (\varphi^{B}_{\alpha})^{-1}
(\{ \xi \in \varphi^B_{\alpha}(U_{\alpha \beta}) \colon \xi_{i_{\alpha}}=0,\ \xi_i>0,\ \text{for}\ i\neq i_{\alpha} \}$. 
The $T^n$-action on $\C^n$ in Example \ref{1.2} preserves 
$\mu_{\C^n}^{-1}(\varphi^{B}_{\alpha}(b))$ which is naturally equivariantly 
diffeomorphic to $T^n/S^1_{i_{\alpha}}$ with the natural $T^n$-action, where $S^1_{i_{\alpha}}$ 
is the circle subgroup 
generated by $e_{i_{\alpha}}$. Under this identification, the restriction 
$(\nu_{\C^n})_{\varphi^{B}_{\alpha}(b)}: \pr_1^{-1}(\varphi^{B}_{\alpha}(b))\to 
\mu_{\C^n}^{-1}(\varphi^{B}_{\alpha}(b))$ of $\nu_{\C^n}$ to the fiber at 
$\varphi^{B}_{\alpha}(b)$ can be thought of as the natural projection from $T^n$ to 
$T^n/S^1_{i_{\alpha}}$. The same argument allows us to identify 
$(\nu_{\C^n})_{\varphi^{B}_{\beta}(b)}: \pr_1^{-1}(\varphi^{B}_{\beta}(b))\to 
\mu_{\C^n}^{-1}(\varphi^{B}_{\beta}(b))$ with the natural projection from 
$T^n$ to $T^n/S^1_{i_{\beta}}$. 
The commutativity $\nu_{\C^n}\circ (\varphi_{\beta}^{B}\times \id_{T^n})
\circ \varphi^T_{\beta \alpha}\circ (\varphi^{B}_{\alpha}\times \id_{T^n})^{-1}=
\varphi^X_{\beta \alpha}\circ \nu_{\C^n}$ implies 
$g_{\beta \alpha}(b)(S^1_{i_{\alpha}})=S^1_{i_{\beta}}$, and since 
the integral lattice of $\Lie (S^1_{i_j})$ is equal to $\Z_{i_j}$ and 
$g_{\beta \alpha}(b)\in SL_n(\Z )\subset \Aut (T^n)$, by taking the derivative 
$g_{\beta \alpha}(b): \Lie (S^1_{i_{\alpha}})\to \Lie (S^1_{i_{\beta}})$, 
$g_{\beta \alpha}(b)$ sends $\Z_{i_{\alpha}}$ to $\Z_{i_{\beta}}$. 
\end{proof}

{\it Proof of Theorem} \ref{thm5.3}. 
Let $U_{\alpha}$ be a locally toric chart with $U_{\alpha}\cap {\mathcal S}^{(n-1)}B\neq \emptyset$. If necessary, by shrinking $U_{\alpha}$, we can assume that 
there exists a unique $i_{\alpha}$ such that 
\[
U_{\alpha}\cap{\mathcal S}^{(n-1)}B=
(\varphi^B_{\alpha})^{-1}(\varphi^B_{\alpha}(U_{\alpha})\cap \{ \xi \in \t^* \colon 
\xi_{i_{\alpha}}=0 \}). 
\]
Let us consider the trivial rank one sub-lattice 
bundle $U_{\alpha}\cap {\mathcal S}^{(n-1)}B\times \Z_{i_{\alpha}}
\to U_{\alpha}\cap {\mathcal S}^{(n-1)}B$ of $\Z^n_P\big|_{U_{\alpha}\cap 
{\mathcal S}^{(n-1)}B}$ 
on $U_{\alpha}\cap {\mathcal S}^{(n-1)}B$, where $\Z_{i_{\alpha}}$ is the rank one 
sub-lattice of $\Z^n$ which is generated by $i_{\alpha}$th fundamental vector $e_{i_{\alpha}}$. 
Lemma \ref{4.1} implies that these sub-lattice bundles on locally toric charts 
which satisfy the above conditions are patched together to the rank one 
sub-lattice bundle $\pi_{{\mathcal L}_X}:{\mathcal L}_X\to {\mathcal S}^{(n-1)}B$ 
of $\Z^n_P\big|_{{\mathcal S}^{(n-1)}B}$ on ${\mathcal S}^{(n-1)}B$. 
From the construction, the primitivity of $\pi_{{\mathcal L}_X}:{\mathcal L}_X\to 
{\mathcal S}^{(n-1)}B$ is obvious. 
\hspace*{9.1cm}
$\Box$\\ 
\begin{ex}
Consider the characteristic bundle of the twisted toric manifold in Example \ref{ex3.7}. Let 
$\pi :\partial \widetilde{B}(=\R \times \{ 0,1\} )\to \partial B(=S^1\times \{ 0,1\} )$ 
be the universal covering space of the boundary $\partial B$ of $B$. 
The restriction of $\pi_{\Z}:\Z^2_P\to B$ to $\partial B$ is the $\Z^2$-bundle 
$\partial \widetilde{B}\times_{\rho }\Z^2$ associated with 
$\pi :\partial \widetilde{B}\to \partial B$ by the homomorphism $\rho:\Z \to SL_2(\Z )$ 
defined in Example \ref{ex3.7}. The primitive rank one sub-lattice $\{ 0\}\times \Z$ of $\Z^2$ is 
preserved by the action $\rho$, and in this case, the characteristic bundle is the associated bundle 
$\pi :\partial \widetilde{B}\times_{\rho }\left( \{ 0\}\times \Z \right) \to \partial B$ 
by the induced $\Z$-action on $\{ 0\}\times \Z$. 
\end{ex}
\begin{ex}
In the case of Example \ref{ex3.9}, the restriction of the associated $\Z^2$-bundle 
$\pi_{\Z }:\Z^2_P\to B$ to the neighborhood $B_2$ of $\partial B$ is obtained from 
the trivial bundle $\pr_1:\overline{B}_2\times \Z^2 \to \overline{B}_2$ by the similar way 
as in the case of the construction of the restriction $T^2_P\big|_{B_2}$, 
that is, by gluing $U_1\times \Z^2$ and $U_2\times \Z^2$ with the diffeomorphism 
$\varphi^{\Z}:U_1\times \Z^2 \to U_2\times \Z^2$ 
\[
\varphi^{\Z}(\xi ,l)=(\varphi^B(\xi ), \rho(\gamma )l). 
\]
The characteristic bundle is obtained by gluing two trivial sub-lattice bundles 
\[
\begin{split}
&\pr_1:\{ \xi \in \overline{B}_2\colon 0<\xi_1,\ \xi_2=0\} \times \left( \{ 0\} \times \Z \right) 
\to \{ \xi \in \overline{B}_2\colon 0<\xi_1,\ \xi_2=0\} , \\
&\pr_1:\{ \xi \in \overline{B}_2\colon \xi_1=0,\ 0<\xi_2\} \times \left( \Z \times \{ 0\} \right) 
\to \{ \xi \in \overline{B}_2\colon \xi_1=0,\ 0<\xi_2\} 
\end{split}
\]
of $\pr_1:\overline{B}_2\times \Z^2 \to \overline{B}_2$ with restrictions 
\[
\begin{split}
&\varphi^{\Z}\big|_{\{ \xi \in U_1\colon \xi_2=0\} \times (\{ 0\} \times \Z )}: 
\{ \xi \in U_1\colon \xi_2=0\} \times \left( \{ 0\} \times \Z \right) \\
&\hspace{7cm}\to \{ \xi \in U_2\colon \xi_1=0\} \times \left( \Z \times \{ 0\} \right) , \\
&\varphi^B\big|_{\{ \xi \in U_1\colon \xi_2=0\}}:\{ \xi \in U_1\colon \xi_2=0\} \to 
\{ \xi \in U_2\colon \xi_1=0\}  
\end{split}
\]
of diffeomorphisms $\varphi^{\Z}$, $\varphi^B$, respectively. 
\end{ex}

Fix a principal $SL_n(\Z)$-bundle $\pi_P:P\to B$ on an $n$-dimensional 
manifold $B$ with corners. 
\begin{defn}
Two twisted toric manifolds $\{ X_1, \nu_1,\mu_1\}$ and $\{ X_2, \nu_2,\mu_2\}$ 
associated with $\pi_P:P\to B$ are {\it topologically isomorphic}, if there exist 
an automorphism $\psi^P$ of $\pi_P:P\to B$ which covers identity map of $B$ 
(then $\psi^P$ induces the automorphism $\psi^T$ of $\pi_T:T^2_P\to B$), and 
a homeomorphism $\psi^X$ from $X_1$ to $X_2$ such that the following diagram 
commutes 
\[
\xymatrix{
T^2_P\ar[rrr]^{\psi^T}\ar[drr]^{\nu_1}\ar[ddr]_{\pi_T} &  
&  & T^2_P\ar[drr]^{\nu_2}\ar[ddr]|\hole ^{\pi_T} 
&  &  \\
  &  & X_1\ar[rrr]^{\psi^X\hspace*{1.4cm}}\ar[dl]^{\mu_1} &  
  &  & X_2\ar[dl]^{\mu_2} \\
  & B \ar@{=}^{\id_B}[rrr] &  &  & B. &  
} 
\]
Note that since the restriction $\nu_i\big|_{\pi_T^{-1}(B\backslash \partial B)}:
\pi_T^{-1}(B\backslash \partial B)\to \mu_i^{-1}(B\backslash \partial B)$ is a 
diffeomorphism for $i=1,2$, so is 
$\psi^X\big|_{\mu_1^{-1}(B\backslash \partial B)}:
\mu_1^{-1}(B\backslash \partial B)\to \mu_2^{-1}(B\backslash \partial B)$.  
\end{defn}

Let $\{ (U_{\alpha}, \varphi^B_{\alpha})\}_{\alpha \in \Gamma}$ be a 
coordinate neighborhood system of $B$. In order to show the 
classification theorem, we assume the following technical conditions 
\begin{enumerate}
\item each coordinate neighborhood $(U_{\alpha}, \varphi^B_{\alpha})$ satisfies the 
condition (i) of Definition \ref{3.1}, 
\item  on a overlap $U_{\alpha \beta}=U_{\alpha}\cap U_{\beta}$ with 
$U_{\alpha \beta}\cap \partial B\neq \emptyset$, 
if $\varphi^B_{\beta \alpha}$ sends 
$\{ \xi \in \varphi_{\alpha}(U_{\alpha \beta})\colon \xi_i=0\}$ to 
$\{ \zeta \in \varphi_{\beta}(U_{\alpha \beta})\colon \zeta_j=0\}$ for some 
$i$, $j$, then $\varphi^B_{\beta \alpha}$ satisfies 
\[
\varphi^B_{\beta \alpha}(\xi)_j=\< e_j, \varphi^B_{\beta \alpha}(\xi )\> =\xi_i 
\] 
on a sufficiently small neighborhood of 
$\{ \xi \in \varphi_{\alpha}(U_{\alpha \beta})\colon \xi_i=0\}$ of 
$\varphi_{\alpha}(U_{\alpha \beta})$. 
\end{enumerate}
In the case where $B$ is a surface, $B$ has such a coordinate neighborhood system. 
In fact, by taking a coordinate neighborhood near $\partial B$ as in Example 
\ref{ex3.9}, we can adopt the composition of a rotation and a parallel transport in 
$\t^*=\R^2$ as a coordinate changing function.  
\begin{thm}\label{classification}
Assume $B$ has a coordinate neighborhood system with the above properties. 
Then by associating the characteristic bundle to a twisted toric manifold, 
the set of topologically isomorphism classes of twisted toric manifolds associated with 
$\pi_P:P\to B$ corresponds one-to-one to the set of equivalent classes of 
primitive rank one sub-lattice bundles on $\CS^{(n-1)}B$ of 
$\Z^n_P\big|_{\CS^{(n-1)}B}$ by the action of the automorphism group of $P$. 
\end{thm}
\begin{proof}
It is easy to see that the map which associates to an isomorphism class of a 
twisted toric manifold the equivalence class of the characteristic bundle is 
well-defined. For two twisted toric manifolds $X_1$ and $X_2$ whose 
characteristic bundles are in the same equivalent class, there exists an 
automorphism $\psi^P$ of $P$ such that induced automorphism $\psi^{\Z}$ 
of $\Z^n_P$ sends the characteristic bundle of $X_1$ to that of $X_2$. 
Then it is easy to see from the construction of the characteristic bundle 
that the automorphism $\psi^T$ of $T^n_P$ induced by $\psi^P$ 
descends to the topologically isomorphism $\psi^X$ from $X_1$ to $X_2$. 
This implies the injectivity. 

Conversely, for each primitive rank one sub-lattice bundle 
$\pi_{\mathcal L}:{\mathcal L}\to {\mathcal S}^{(n-1)}B$ 
of $\Z^n_P\big|_{{\mathcal S}^{(n-1)}B}$ on ${\mathcal S}^{(n-1)}B$, 
we can construct a twisted toric manifold whose characteristic bundle is 
equal to ${\mathcal L}$ in the following way. 
Let $(U, \varphi^{B})$ be a coordinate neighborhood 
of $B$ which satisfies the condition (i) of Definition \ref{3.1}. We set 
\[
k=\# \{ i \colon \exists \xi \in \varphi^{B}(U)\subset \t^*_{\ge 0}\ 
s.t.\ \xi_i=0 \} . 
\]
If $k=0$, just define $\mu^{-1}(U)=\pi_T^{-1}(U)$. If $k>0$, let 
\[
\{ i_1,\ldots ,i_k \} =\{ i \colon \exists \xi \in \varphi^{B}(U)\ s.t.\ \xi_i=0 \} 
\]
and we denote by $(U\cap {\mathcal S}^{(n-1)}B)_a$ the connected 
component of $U\cap {\mathcal S}^{(n-1)}B$ which satisfies 
\[
\varphi^{B} ((U\cap {\mathcal S}^{(n-1)}B)_a)
=\{ \xi \in \varphi^{B} (U) \colon \xi_{i_a}=0,\ \xi_i>0\ \text{for}\ i\neq i_a\}  
\] 
for $a=1, \ldots , k$. Since $\pi_{\mathcal L}:{\mathcal L}\to {\mathcal S}^{(n-1)}B$ 
is primitive and $U$ is contractible, there is a local trivialization 
$\varphi^P:\pi_P^{-1}(U)\to U\times SL_n(\Z )$ 
(hence $\varphi^P$ induces $\varphi^{\Z}:\pi_{\Z}^{-1}(U)\to U\times \Z^n$)
and a primitive tuple $\{ L_1, \ldots , L_k \}$ of rank one sub-lattices of $\Z^n$ 
such that $\varphi^{\Z}$ trivialize 
$\pi_{{\mathcal L}}^{-1}((U\cap {\mathcal S}^{(n-1)}B)_a)\cong 
(U\cap {\mathcal S}^{(n-1)}B)_a\times L_a$. If necessary, 
by applying an element of $SL_n(\Z )$ to $\{ L_1, \ldots , L_k \}$, 
we may assume that $L_a$ is equal to the rank one sub-lattice $\Z_{i_a}$ 
in $\Z^n$ spanned by the $i_a$th fundamental vector $e_{i_a}$ 
for $a=1, \ldots ,k$. This is possible because $\{ L_1, \ldots , L_k \}$ is primitive. 
Let $\widetilde{\varphi}^T:\pi_T^{-1}(U)\cong \varphi^{B}(U)\times T^n$ 
be the composition of the trivialization 
$\varphi^T:\pi_T^{-1}(U)\cong U\times T^n$ which is induced by $\varphi^P$ and 
$\varphi^{B}\times \id_{T^n}$. Consider the composition of $\widetilde{\varphi}^T$ 
and the map $\nu_{\C^n}:\varphi^{B}(U)\times T^n \to 
\mu_{\C^n}^{-1}(\varphi^{B}(U))$ which is defined in (\ref{eq3.1}). 
On $U\backslash \partial B$, this gives the orientation preserving 
diffeomorphism from $\pi_T^{-1}(U\backslash \partial B)$ to 
$\mu_{\C^n}^{-1}(\varphi^{B}(U\backslash \partial B)$. We define 
$\mu^{-1}(U)$ to be the space obtained by gluing 
$\pi_T^{-1}(U\backslash \partial B)$ and $\mu_{\C^n}^{-1}(\varphi^{B}(U))$ 
by this map. 

Let $\{ (U_{\alpha}, \varphi^B_{\alpha})\}_{\alpha \in \Gamma}$ be a 
coordinate neighborhood system of $B$ which satisfies the assumption before 
the theorem. We apply 
the above construction to each $(U_{\alpha}, \varphi^B_{\alpha})$.  
For each overlap with $U_{\alpha}\cap U_{\beta}\cap \partial B\neq \emptyset$, 
it is easy to check from the second condition in the assumption,   
Lemma \ref{B2}, and Lemma \ref{a4} that the transition map 
$\widetilde{\varphi}^T_{\beta \alpha}
:\varphi^{B}_{\alpha}(U_{\alpha \beta})\times T^n\cong 
\varphi^{B}_{\beta}(U_{\alpha \beta})\times T^n$ defined in Section 4.2 
descends to the orientation preserving diffeomorphism from 
$\mu_{\C^n}^{-1}(\varphi^{B}_{\alpha}(U_{\alpha \beta}))$ to 
$\mu_{\C^n}^{-1}(\varphi^{B}_{\beta}(U_{\alpha \beta}))$. 
Then we can glue $\mu^{-1}(U_{\alpha})$ for all $U_{\alpha}$ to obtain the 
twisted toric manifold $X$ whose characteristic bundle is equal to ${\mathcal L}$. 
This implies the surjectivity. 
\end{proof}

\section{Topology}
Let $\pi_P:P\to B$ be a principal $SL_n(\Z)$-bundle on a $n$-dimensional 
manifold $B$, $\pi_T:T^n_P\to B$ the $T^n$-bundle associated with $P$ by 
the natural action of $SL_n(\Z )$ on $T^n$, $X$ a $2n$-dimensional 
twisted toric manifold associated with $P$.  
\subsection{Fundamental groups}
In this section, we shall investigate the fundamental group of X. Let 
$b\in B$ be the base point which is located in the interior of $B$. 
Fix the base point $y_0\in \pi_T^{-1}(b)$ and set $x_0=\nu (y_0)\in X$.  
Since $T^n$ is identified with $\R^n /\Z^n$, $\pi_T:T^n_P\to B$ has 
the zero-section $s_T$. 
Then the homotopy exact sequence for $T^n_P$ 
splits into the short exact sequence, and $\pi_1(T^n_P, y_0)$ is isomorphic 
to the semi-direct product $\pi_1(\pi_T^{-1}(b), y_0)\rtimes \pi_1(B ,b)$ 
of $\pi_1(\pi_T^{-1}(b),y_0)$ and $\pi_1(B ,b)$. 
The section $s_T$ also defines the section $s_X$ of 
$\mu :X\to B$ by $s_X=\nu \circ s_T$ 
which also gives the identification of $\pi_1(X, x_0)$ with the semi-direct 
product $\ker \mu_* \rtimes \pi_1(B ,b)$.   

Let us consider the following commutative diagram of homomorphisms 
\[
\xymatrix{
\{ 1\} \ar[r] & \pi_1(\pi_T^{-1}(b),y_0)\ar[r]^{\iota_*} \ar[d]^{\kappa} & \pi_1(T^n_P, y_0) \ar[r]^{{\pi_T}_*} \ar[d]^{\nu_*} & 
\pi_1(B ,b) \ar[r] \ar@{=}[d]& \{ 1\} \ \ \text{(exact)}\\
\{ 1\} \ar[r] & \ker \mu_*\ar[r] & \pi_1(X, x_0) \ar[r]^{\mu_*} & \pi_1(B ,b) \ar[r] & \{ 1\} \ \ 
\text{(exact),} 
}
\]
where $\iota :\pi_T^{-1}(b)\hookrightarrow T^n_P$ is the inclusion map and 
$\kappa :\pi_1(\pi_T^{-1}(b),y_0)\to \ker \mu_*$ is defined by 
$\kappa =\nu_*\circ \iota_*$. 
\begin{lem}
The map $\kappa : \pi_1(\pi_T^{-1}(b), y_0)\to \ker \mu_*$ is surjective. 
\end{lem}
Note that the surjectivity of the map $\kappa$ is equivalent to that of 
$\nu_*$ since ${s_X}_*=\nu_*\circ {s_T}_*$.
\begin{proof}
Since $b$ is in the interior of $B$, the composition map 
$\nu \circ \iota :\pi_T^{-1}(b)\to \mu^{-1}(b)$ is a 
diffeomorpshims which sends $y_0$ to $x_0$.  
Then it is sufficient to show that every element of $\ker \mu_*$ is 
represented by the loop in $\mu^{-1}(b)$. Let $\alpha \in \ker \mu_*$ 
and take its representative $a':I \to X$ with $a'(0)=a'(1)=x_0$. 
This means that the map $\mu \circ a':I\to B$ is homotopic to the constant 
map $b$. Then if necessary, by replacing a representative of $\alpha$, 
we can assume that there exists a contractible locally toric chart $U$ 
located in the interior of $B$ such that the image of $a'$ is included in 
$\mu^{-1}(U)$. Since $U$ is a locally toric chart which is included in the 
interior of $B$, $\mu^{-1}(U)$ is diffeomorphic to $U\times T^n$, and the fact 
that $U$ is contractible implies that $\mu^{-1}(b)$ is a deformation retract 
of $\mu^{-1}(U)$. We take a deformation map 
$h:I\times \mu^{-1}(U)\to \mu^{-1}(U)$ 
which satisfies $h(0, \cdot)=id_{\mu^{-1}(U)}$ and 
$h(1, \cdot): \mu^{-1}(U)\to \mu^{-1}(b)$. Then the map $h(s, a'(t))$ is 
the homotopy connecting $a'$ with the loop $a=h(1, a'(t))$ in $\mu^{-1}(b)$. 
\end{proof}

Since $SL_n(\Z )$ is discrete, $P$ is a flat $SL_n(\Z)$-bundle. Then, 
$T^n_P$ also has the flat connection which is induced from that of $P$.
Let $\gamma :I\to B $ be a path of $B$ which starts from the base point 
$b$, and ${\rm Hol}^T_{\gamma}:\pi_T^{-1}(b)\to \pi_T^{-1}(\gamma (1))$ 
the parallel transport of $T^n_P$ with respect to the induced flat connection 
along $\gamma$. For $\gamma :I\to B$, define the subset $K_{\gamma}$ of 
$\pi_1(\pi_T^{-1}(b),y_0)$ by 
\[
 K_{\gamma}=\{ \alpha \in \pi_1(\pi_T^{-1}(b) ,y_0) \colon (\nu \circ {\rm
 Hol}_{\gamma}^T)_*(\alpha )=1  \} .
\]
We denote by $K$ the subgroup of $\pi_1(\pi_T^{-1}(b),y_0)$ which is 
generated by $\displaystyle \cup_{\gamma}K_{\gamma}$ where $\gamma$ runs 
over all paths of $B$ which start at $b$. Then the following lemma is 
obvious. 
\begin{lem}
The subgroup $K$ is included in the kernel of $\kappa$. 
\end{lem}
\begin{thm}\label{fundamental}
If $B$ has at least one corner point, then $\mu :X\to B$ induces the 
isomorphism $\mu_*:\pi_1(X,x_0)\cong \pi_1(B,b)$. 
\end{thm}
\begin{proof}
Let us show $K=\pi_1(\pi_T^{-1}(b),y_0)$. Let $\gamma :I\to B$ be a 
path which starts at $b$ to the corner point $b'$ of $B$. Since 
$\mu^{-1}(b')$ consists of only one point, the composition map 
$\nu \circ {\rm Hol}^T_{\gamma}:\pi_T^{-1}(b)\to \mu^{-1}(b')$ sends 
every loop of $\pi_T^{-1}(b)$ to the constant map. This implies 
$K=\pi_1(\pi_T^{-1}(b),y_0)$. 
\end{proof}
\begin{rem}
From facts that the homomorphism $\mu_*:\pi_1(X, x_0)\to \pi_1(B ,b)$ of 
fundamental groups induced by $\mu :X\to B$ is surjective and that complete 
non-singular toric varieties in the original algebro-geometric sense are simply 
connected \cite{F}, we can see that twisted toric manifolds associated with 
principal $SL_n(\Z )$-bundles on non-simply connected manifolds, such as 
those in Example \ref{ex3.7} and \ref{ex3.9}, are not toric varieties. 
In particular, these are not symplectic toric manifolds.  
\end{rem}

\subsection{Cohomology groups}
Since $\mu :X\to B$ has a section, as described in the last subsection, 
the induced homomorphism $\mu^*:H^1(B;\Z )\to H^1(X;\Z )$ is injective.   
In particular, if $B$ has at least one corner point, $\mu^*$ is isomorphism. 

In this section, we shall give the method to calculate not only the first cohomology 
but also the full cohomology group of a twisted toric manifold. For topological tools 
which we use in this section, see \cite{H}. 
Let $X$ be a $2n$-dimensional twisted toric manifold associated with a 
principal $SL_n(\Z)$-bundle on a $n$-dimensional manifold $B$, and 
$T^n_P$ the $T^n$-bundle associated with $P$ by the natural 
action of $SL_n(\Z )$ on $T^n$. Note that since $SL_n(\Z )$ is discrete, 
$P$ is flat, hence this induces the flat connection to $T^n_P$. 

Assume $B$ is equipped with a CW complex structure. For each 
$p$-dimensional cell $e^{(p)}$ (we often say simply $p$-cell), let $c$ be the barycenter 
of $e^{(p)}$, the map $\varphi:\overline{D^p}\to B$ denotes 
the characteristic map of $e^{(p)}$ from the $p$-dimensional closed ball 
$\overline{D^p}$ to $B$. 
Define the map $\widetilde{\varphi}^T:
\overline{D^p}\times \pi_T^{-1}(c)\to T^n_P$ by 
\[
\widetilde{\varphi}^T(d, \theta )=
{\rm Hol}_{\varphi \circ \gamma }^T(\theta )
\]
for $(d, \theta )\in \overline{D^p}\times \pi_T^{-1}(c)$, where  
$\gamma :[0, 1]\to \overline{D^p}$ is a path with $\gamma (0)=c$ 
and $\gamma (1)=d$ and 
${\rm Hol}_{\varphi \circ \gamma }^T:(T^n_P)_c\to 
(T^n_P)_d$ is the parallel transport of $T^n_P$ with respect to the induced connection 
along $\varphi \circ \gamma :[0, 1]\to B$. Note that 
${\rm Hol}_{\varphi \circ \gamma }^T$ does not depend on the 
choice of $\gamma$, since $\overline{D^p}$ is contractible.   
We shall assume that the cell decomposition of $B$ satisfies the following 
conditions 
\begin{enumerate}
\item the restriction of $\mu :X\to B$ to each $p$-dimensional cell $e^{(p)}$ 
is a (trivial) torus bundle, 
\item for each $p$-dimensional cell $e^{(p)}$, the map $\widetilde{\varphi}^T$ 
induces the map $\widetilde{\varphi}^X:\overline{D^p}\times \mu^{-1}(c)\to X$ 
such that the following diagram commutes 
\[
\xymatrix{
\overline{D^p}\times \pi_T^{-1}(c)
\ar[rrrr]^{\id_{\overline{D^p}}\times \nu}\ar[drr]_{\pr_1}
\ar[ddr]_{\widetilde{\varphi}^T} & &   & & 
\overline{D^p}\times \mu^{-1}(c)
\ar[dll]^{\pr_1}\ar[ddr]^{\widetilde{\varphi}^X} & \\
                                          & & \overline{D^p}
										  \ar[ddr]^{\varphi }|\hole & & & \\
    & T^n_P\ar[rrrr]^{\nu}\ar[drr]_{\pi_T} & & & &  X\ar[dll]^{\mu} \\
    &                              & & B. & &                     
}
\] 
\end{enumerate}
\begin{rem}\label{rem6.1}
This assumption is achieved when each $p$-cell $e^{(p)}$ is included in a 
$k$-dimensional strata ${\mathcal S}^{(k)}B$ of $B$. 
%
Since the characteristic map $\varphi$ sends the interior $D^p$ 
of $\overline{D^p}$ homeomorphically to $e^{(p)}$, $k$ is necessarily equal or 
greater than $p$.  
\end{rem}

In the rest of this section, {\it we shall assume the condition in Remark} \ref{rem6.1}.  
Let $B^{(p)}$ be the $p$-skeleton of $B$, and 
$(T^n_P)^{(p)}=\pi_T^{-1}(B^{(p)} )$, $X^{(p)}=\mu^{-1}(B^{(p)})$ 
its inverse images by $\pi_T:T^n_P\to B$, $\mu :X\to B$, respectively.   
We consider the spectral sequence $\{ (E_X)^{p,q}_r, d^X_r \}$ with respect to the 
filtration 
\[
S^*(X;\Z )\supset S^*(X, X^{(0)}; \Z )\supset \cdots 
\supset S^*(X, X^{(n)};\Z )=0 \\
\]
of the singular cochain complex with coefficient $\Z$. $\{ (E_X)^{p,q}_r, d^X_r \}$ 
is called the {\it cohomology Leray spectral sequence} of the map $\mu :X\to B$.  
Let $c^{(p)}_{\lambda}$ be the barycenter of the $p$-cell $e^{(p)}_{\lambda}$. 
As mentioned before, we can identify the fibers $\pi_T^{-1}(c^{(p)}_{\lambda})$ 
and $\mu^{-1}(c^{(p)}_{\lambda})$ with tori. Then 
$\nu_{c^{(p)}_{\lambda}}:\pi_T^{-1}(c^{(p)}_{\lambda})\to 
\mu^{-1}(c^{(p)}_{\lambda})$ can be thought of as a surjective homomorphism
between them, and we have the exact sequence of tori
\[
0\to \ker \nu_{c^{(p)}_{\lambda}} \hookrightarrow \pi_T^{-1}(c^{(p)}_{\lambda}) 
\stackrel{\nu_{c^{(p)}_{\lambda}}}{\to} \mu^{-1}(c^{(p)}_{\lambda})\to 0.  
\]
Since any exact sequence of tori splits, the map $\nu_{c^{(p)}_{\lambda}}$ induces 
the injective homomorphism 
$\nu_{c^{(p)}_{\lambda}}^*:H^*(\mu^{-1}(c^{(p)}_{\lambda}) ;\Z )\hookrightarrow  
H^*(\pi_T^{-1}(c^{(p)}_{\lambda}) ;\Z )$ which enable us to identify 
$H^q(\mu^{-1}(c^{(p)}_{\lambda}) ;\Z )$ with its image in 
$H^q(\pi_T^{-1}(c^{(p)}_{\lambda}) ;\Z )$ by $\nu_{c^{(p)}_{\lambda}}^*$. 

Let $\left( C^p(B ; {\mathcal H}_T^q), \delta \right)$ be the cochain complex 
of the CW complex $B$ with the Serre local system ${\mathcal H}_T^q$ of the $q$th 
cohomology with $\Z$-coefficient for the torus bundle $\pi_T: T^n_P\to B$. 
We denote by $C^p(B ; {\mathcal H}_X^q)$ the subset of 
$C^p(B ; {\mathcal H}_T^q)$ whose cochain takes a value in the image 
$\nu_{c^{(p)}_{\lambda}}^*\left( H^q(\mu^{-1}(c^{(p)}_{\lambda});\Z )\right)$ 
of $q$th cohomology $H^q(\mu^{-1}(c^{(p)}_{\lambda});\Z )$ by 
$\nu_{c^{(p)}_{\lambda}}^*$ for each $p$-cell $e^{(p)}_{\lambda}$. 
\begin{thm}\label{thm6.2}
$C^p(B ; {\mathcal H}_X^q)$ is preserved by the differential $\delta$ of 
$C^p(B ; {\mathcal H}_T^q)$, so that $C^p(B ; {\mathcal H}_X^q)$ is a 
sub-complex of $\left( C^p(B ; {\mathcal H}_T^q), \delta \right)$. We denote 
its cohomology by $H^p(B ;{\mathcal H}^q_X)$. Then we have the isomorphisms 
\begin{align}
(E_X)_1^{p, q}&\cong C^p(B ; {\mathcal H}^q_X), \ \ \ 
(E_X)_2^{p, q}\cong H^p(B ; {\mathcal H}^q_X), \nonumber \\ 
(E_X)_{\infty}^{p,q} &=F^pH^{p+q}(X; \Z )/F^{p+1}H^{p+q}(X;\Z ), \nonumber 
\end{align}
where $F^lH^k(X; \Z )$ is the image of the map $H^k(X, X^{(l-1)};\Z )\to H^k(X;\Z )$. 
\end{thm}
\begin{proof}
Let $\{ (E_T)^{p,q}_r, d^T_r\}$ be the cohomology Serre spectral sequence of the 
torus bundle $\pi_T:T^n_P\to B$, that is, the spectral sequence with respect to 
the filtration 
\[
S^*(T^n_P;\Z )\supset 
S^*(T^n_P, (T^n_P)^{(0)}; \Z )\supset \cdots \supset 
S^*(T^n_P, (T^n_P)^{(n)};\Z )=0 
\] 
of the singular cochain complex with coefficient $\Z$. 
By the excision isomorphism, the above assumption, and the K\"unneth 
formula, for $E_1$-terms, we have the isomorphisms 
\begin{align} 
(E_T)_1^{p,q}&=H^{p+q}((T^n_P)^{(p)}, (T^n_P)^{(p-1)}; \Z ) \nonumber \\
&\cong \sum_{\lambda}H^{p+q}(\pi_T^{-1}(\overline{e^{(p)}_{\lambda}}), 
\pi_T^{-1}(\overline{e^{(p)}_{\lambda}}-e^{(p)}_{\lambda});\Z ) \nonumber \\
&\cong \sum_{\lambda}H^{p+q}((\overline{D^p}_{\lambda}, 
\partial \overline{D^p}_{\lambda})\times \pi_T^{-1}(c^{(p)}_{\lambda}); \Z )\nonumber \\
&\cong \sum_{\lambda}H^p(\overline{D^p}_{\lambda}, 
\partial \overline{D^p}_{\lambda} ;\Z )\otimes H^q(\pi_T^{-1}(c^{(p)}_{\lambda});\Z )
, \nonumber 
\end{align}
\begin{align}
(E_X)_1^{p,q}&=H^{p+q}(X^{(p)}, X^{(p-1)}; \Z ) \nonumber \\
&\cong \sum_{\lambda}H^{p+q}(\mu^{-1}(\overline{e^{(p)}_{\lambda}}), 
\mu^{-1}(\overline{e^{(p)}_{\lambda}}-e^{(p)}_{\lambda});\Z ) \nonumber \\
&\cong \sum_{\lambda}H^{p+q}((\overline{D^p}_{\lambda}, 
\partial \overline{D^p}_{\lambda})\times \mu^{-1}(c^{(p)}_{\lambda}); \Z )\nonumber \\
&\cong \sum_{\lambda}H^p(\overline{D^p}_{\lambda}, 
\partial \overline{D^p}_{\lambda} ;\Z )\otimes H^q(\mu^{-1}(c^{(p)}_{\lambda});\Z )
, \nonumber 
\end{align}
where the sum runs over all $p$-dimensional cells $e^{(p)}_{\lambda}$. 
On the other hand, by the injectivity of the homomorphism 
$\nu_{c^{(p)}_{\lambda}}^*:H^*(\mu^{-1}(c^{(p)}_{\lambda}) ;\Z )\hookrightarrow  
H^*(\pi_T^{-1}(c^{(p)}_{\lambda}) ;\Z )$ and the assumption (ii), the map 
$\nu:(T^n_P)^{(p)}\to X^{(p)}$ induces the natural injection 
$\nu^*: (E_X)^{p,q}_1\hookrightarrow (E_T)^{p,q}_1$ such that the following diagram 
commutes
\[
\xymatrix@C=0pt{
(E_X)^{p,q}_1\cong \ar@{^{(}->}[d]_{\nu^*}\ar@{}[dr]|{\circlearrowright} & 
\sum_{\lambda}H^p(\overline{D^p}_{\lambda}, 
\partial \overline{D^p}_{\lambda} ;\Z )\otimes H^q(\mu^{-1}(c^{(p)}_{\lambda});\Z )
\ar@{^{(}->}[d]^{\sum_{\lambda}\id^*_{\overline{D^p}} \otimes \nu^*_{c^{(p)}_{\lambda}}} \\
(E_T)^{p,q}_1\cong & 
\sum_{\lambda}H^p(\overline{D^p}_{\lambda}, 
\partial \overline{D^p}_{\lambda} ;\Z )\otimes H^q(\pi_T^{-1}(c^{(p)}_{\lambda});\Z ) . 
}
\]
Moreover it is well known that the $E_1$-term $\{ (E_T)^{p,q}_1, d^T_1\}$ of the 
Serre spectral sequence is isomorphic to the CW complex 
$\left( C^p(B ; {\mathcal H}_T^q), \delta \right)$ with the Serre local system 
${\mathcal H}_T^q$ for the torus bundle $\pi_T: T^n_P\to B$. This fact and 
the naturality of the maps in the spectral sequences prove Theorem \ref{thm6.2}. 
\end{proof}
\begin{rem}
(1) For $q=0$, it is easy to see that $(E_X)_2^{p,0}\cong H^p(B; {\mathcal H}_X^q)\cong H^p(B;\Z )$. 
Moreover $(E_X)_1^{p,q}=0$, if $q$ or $p$ is greater than half the dimension of $X$. \\
(2) If $n=2$ and $\partial B \neq \emptyset$, we can take a cell decomposition of $B$ so that 
all zero cells are included in $\partial B$. In this case, the Leray spectral sequence 
$\{ (E_X)^{p,q}_r, d^X_r \}$ degenerates at $E_2$-term. In fact, $\partial B\neq \emptyset$ 
implies $(E_X)_2^{2, 0}\cong H^2(B; {\mathcal H}_X^0)\cong H^2(B; \Z )=0$, and since 
$e^{(0)}_{\lambda}\in \partial B$, the fiber $\mu^{-1}(e^{(0)}_{\lambda})$ 
of $\mu$ on $e^{(0)}_{\lambda}$ is diffeomorphic to the torus whose dimension is equal or 
less than one. Then $(E_X)_2^{0,2}\cong (E_X)_1^{0,2}\cong C^0(B; {\mathcal H}_X^2)=0$. 
\end{rem}
\begin{cor}
The Euler characteristic $\chi (X)$ is equal to the number of the $0$-dimensional strata 
${\mathcal S}^{(0)}B$ of $B$. 
\end{cor}
\begin{proof}
Let us consider the rational coefficient cohomology Leray spectral sequence 
$\{ (E_X)^{p,q}_r, d^X_r \}$ of the map $\mu :X\to B$.  
Define 
\[
\chi ((E_X)_r)=\sum_{p,q}(-1)^{p+q}\dim_{\Q}(E_X)^{p,q}_r. 
\]
Since $(E_X)^{p,q}_1=C^p(B ;{\mathcal H}^q_X)$, 
\begin{align}
\chi ((E_X)_1)&=\sum_{p,q}(-1)^{p+q}\dim_{\Q}(E_X)^{p,q}_1 \nonumber \\
&=\sum_{p,q}(-1)^{p+q}\dim_{\Q}C^p(B ;{\mathcal H}^q_X) \nonumber \\
&=\sum_{p,q}(-1)^{p+q}\sum_{\lambda}\dim_{\Q}H^q(\mu^{-1}(c^{(p)}_{\lambda}) ;\Q ) 
\nonumber \\
&=\sum_{p}\sum_{\lambda}(-1)^p\chi (\mu^{-1}(c^{(p)}_{\lambda})) , \label{euler}
\end{align}
where the summation $\sum_{\lambda}$ runs over all $p$-cells. 
From the construction of the twisted toric manifold, we have 
\[
\chi (\mu^{-1}(c^{(p)}_{\lambda}))=
\begin{cases}
1 & \ c^{(p)}_{\lambda}\in {\mathcal S}^{(0)}B \\
0 & \ \text{otherwise. }
\end{cases} 
\]
By the assumption of the cell decomposition of $B$, 
$c^{(p)}_{\lambda}\in {\mathcal S}^{(0)}B$ if and only if $c^{(p)}_{\lambda}$ is the 
barycenter of the $0$-cell in ${\mathcal S}^{(0)}B$, in particular, $p=0$. 
Then (\ref{euler}) is equal to the number of ${\mathcal S}^{(0)}B$. 
On the other hand, it is easy to see 
\[
\chi ((E_X)_r)=\chi ((E_X)_1) 
\]
for any $r$, and since $B$ is compact and all fibers of $\mu$ have finitely 
generated cohomology groups, we can obtain  
\[
(E_X)_{\infty}=(E_X)_r
\]
for sufficiently large $r$. Moreover, from 
$(E_X)_{\infty}^{p,q}=F^pH^{p+q}(X; \Q )/F^{p+1}H^{p+q}(X;\Q )$, 
we can easily check that $\chi (X)=\chi ((E_X)_{\infty})$. This proves the corollary. 
\end{proof}

In the rest of this subsection, we shall calculate the cohomology groups for 
some examples. 
\begin{ex}[Example \ref{ex3.7}]\label{ex6.4}
Let us calculate the cohomology group for Example \ref{ex3.7}. 
Let $Q=[0,1]\times [0,1]$ be the square in $\R^2$. 
In this case, $B$ is a cylinder, so $B$ is obtained from $Q$ 
by identifying each point $(0,\xi_2)$ with $(1,\xi_2)$ in $Q$. Then we have 
the natural map from $Q$ to $B$ which is denoted by 
$\varphi :Q\to B$. $Q$ gives a cell decomposition 
of $B$ as in Figure \ref{fig5}.   
\begin{figure}[hbtp]
\begin{center}
\input{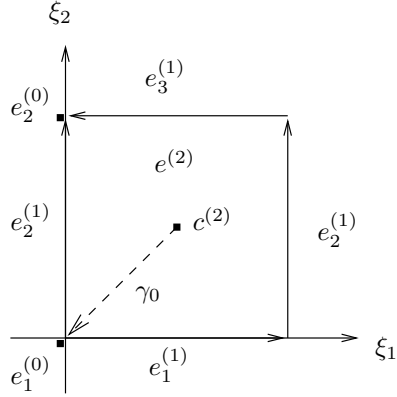}
\caption{$Q$ and cell decomposition}
\label{fig5}
\end{center}
\end{figure}
The pull-back of the triangle of commutative maps $\pi_T:T^2_P\to B$, 
$\nu :T^2_P\to X$, and $\mu :X\to B$ by $\varphi$ are naturally identified with the 
restriction of the triangle (\ref{triangle}) in Example \ref{ex3.7} to $Q$
\[
\xymatrix{
Q\times T^2\ar[rrr]_{\ \ \ \ \ \widetilde{\varphi}^T}\ar[drr]^{\overline{\nu}|_{Q\times T^2}}
\ar[ddr]^{\pr_1} &        
&     & T^2_P\ar[drr]^{\nu}\ar[ddr]|\hole ^{\pi_T} 
&        &     \\
      &        & \overline{\mu}^{-1}(Q)\ar[rrr]_{\widetilde{\varphi}^X\hspace*{5mm}}
\ar[dl]^{\overline{\mu}|_{\overline{\mu}^{-1}(Q)}} &       
      &        & X\ar[dl]^{\mu} \\
      & Q\ar[rrr]^{\varphi} &     &       & B , &     
}
\] 
where $\widetilde{\varphi}^T:Q\times T^2\to T^2_P$ and 
$\widetilde{\varphi}^X:\overline{\mu}^{-1}(Q)\to X$ denote the induced fiberwisely 
diffeomorphisms by pull-back. Under the identifications of fibers by $\widetilde{\varphi}^T$ 
and $\widetilde{\varphi}^X$, the fibers $\mu^{-1}(c^{(p)}_{\lambda})$ on all cells 
$e^{(p)}_{\lambda}$ except for $e^{(1)}_2$ and $e^{(2)}$ are diffeomorphic to 
$T^2/(0\times S^1)$ and the map 
$\nu_{c^{(p)}_{\lambda}}=\nu\big|_{\pi_T^{-1}(c^{(p)}_{\lambda})}:\pi_T^{-1}(c^{(p)}_{\lambda})\to 
\mu^{-1}(c^{(p)}_{\lambda})$ can be identified with the natural projection 
$T^2\to T^2/(0\times S^1)$. Then the image of the induced injection 
$\nu_{c^{(p)}_{\lambda}}^*:H^q(\mu^{-1}(c^{(p)}_{\lambda}); \Z )\to 
H^q(\pi_T^{-1}(c^{(p)}_{\lambda}); \Z )$ can be given by 
\begin{equation}\label{eq6.1}
\nu_{c^{(p)}_{\lambda}}^*\left( H^q(\mu^{-1}(c^{(p)}_{\lambda}); \Z )\right) =
\begin{cases}
\Z & q=0 \\
\Z \oplus 0 & q=1 \\
0 & \text{otherwise. }
\end{cases}
\end{equation}
On the cell $c^{(p)}_{\lambda}=e^{(1)}_2$, $e^{(2)}$, the fiber 
$\mu^{-1}(c^{(p)}_{\lambda})$ is naturally diffeomorphic to 
$\pi_T^{-1}(c^{(p)}_{\lambda})\cong T^2$ by $\nu_{c^{(p)}_{\lambda}}$. 


For $q=0$, the cohomology $H^p(B ; {\mathcal H}_X^0)$ of 
$\left( C^p(B ; {\mathcal H}_X^0), \delta \right)$ is naturally identified with 
the cohomology $H^p(B ; \Z )$ with $\Z$-coefficient. 

For $q=1$, the degree $p$ cochain $u\in C^p(B ; {\mathcal H}_X^1)$ takes a value 
as follows
\[
u(e^{(p)}_{\lambda})\in 
\nu_{c^{(p)}_{\lambda}}^*\left( H^1(\mu^{-1}(c^{(p)}_{\lambda});\Z )\right) =
\begin{cases}
\Z\oplus 0 & p=0,\ \text{or}\ p=1\ \text{and}\ \lambda =1, 3 \\
\Z\oplus \Z & p=1\ \text{and}\ \lambda =2,\ \text{or}\ p=2 \\
0 & \text{otherwise. }
\end{cases}
\]
All differentials $\delta^p:C^p(B;{\mathcal H}_X^1)\to C^{p+1}(B;{\mathcal H}_X^1)$ 
vanish except for $p=0, 1$, and for $p=0$, 
$\delta^0:C^0(B;{\mathcal H}_X^1)\to C^1(B;{\mathcal H}_X^1)$ is given as follows 
\begin{align}
\left( \delta^0u\right) (e^{(1)}_1)&= ^t\rho (1)^{-1}u(e^{(0)}_1)- u(e^{(0)}_1), \nonumber \\
\left( \delta^0u\right) (e^{(1)}_2)&= u(e^{(0)}_2)- u(e^{(0)}_1), \nonumber \\
\left( \delta^0u\right) (e^{(1)}_3)&= ^t\rho (-1)^{-1}u(e^{(0)}_2)- u(e^{(0)}_2). \nonumber  
\end{align}
Fix the path $\gamma_0$ which starts from $c^{(2)}$ to $c^{(0)}_1$. For each 
$c^{(1)}_{\lambda}$, let $\gamma_{\lambda}$ be the path which connects $c^{(0)}_1$ 
and $c^{(1)}_{\lambda}$ counterclockwisely along the boundary of $Q$. 
We identify each fiber $\pi_T^{-1}(c^{(1)}_{\lambda})$ with $\pi_T^{-1}(c^{(2)})$ by 
the parallel transport along the composition of the paths $\gamma_0$ and 
$\gamma_{\lambda}$ with respect to the connection induced from that of $P$. Then 
$\delta^1:C^1(B;{\mathcal H}_X^1)\to C^2(B;{\mathcal H}_X^1)$ is given by
\[
\left( \delta^1u\right) (e^{(2)})= u(e^{(1)}_1)+ ^t\rho (1)^{-1}u(e^{(1)}_2)
+ ^t\rho (1)^{-1}u(e^{(1)}_3) - ^t{\rho (1)^{-1}} ^t\rho (-1)^{-1}u(e^{(1)}_2). 
\]
Then the cohomology groups are obtained by 
\[
H^p(B ; {\mathcal H}_X^1)=
\begin{cases}
\Z /2\Z & p=1, 2 \\
0 & \text{otherwise. }
\end{cases}
\]

For $q=2$, 
the degree $p$ cochain $u\in C^p(B ; {\mathcal H}_X^2)$ takes a value as follows
\[
u(e^{(p)}_{\lambda})\in 
\nu_{c^{(p)}_{\lambda}}^*\left( H^1(\mu^{-1}(c^{(p)}_{\lambda});\Z )\right) =
\begin{cases}
\Z & p=1\ \text{and}\ \lambda =2,\ \text{or}\ p=2 \\
0 & \text{otherwise.}   
\end{cases}
\]
In this case, all differentials 
$\delta^p:C^p(B ; {\mathcal H}_X^2)\to C^{p+1}(B ; {\mathcal H}_X^2)$ vanish. 
It is clear except for $p=1$. For $p=1$, 
since the holonomy $\rho (-1)$ along $e^{(1)}_1$ 
(resp. $\rho (1)$ along $e^{(1)}_3$) induces the identity of $H^2(\pi_T^{-1}(c^{(0)}_1);\Z )$ (resp. $H^2(\pi_T^{-1}(c^{(0)}_2);\Z )$), the differential $\delta^1:C^1\to C^2$ is given by  
\[
\delta^1u(e^{(2)})=u(e^{(1)}_1)+u(e^{(1)}_2)+u(e^{(1)}_3)- u(e^{(1)}_2) 
  = u(e^{(1)}_2)-u(e^{(1)}_2)=0. 
\]
Then the cohomology groups are obtained by 
\[
H^p(B ; {\mathcal H}_X^2)=
\begin{cases}
\Z & p=1, 2 \\
0 & \text{otherwise. }
\end{cases}
\]
The table for the $E_2$-terms is in Figure \ref{fig6}. 
\begin{figure}[hbtp]
\begin{center}
\input{fig6.pstex_t}
\caption{the table of $(E_X^{p,q})_2$-terms for Example \ref{ex3.7}}
\label{fig6}
\end{center}
\end{figure}
In particular, the Leray spectral sequence is 
degenerate at $E^2$-term in this case, and the cohomology groups of $X$ are given by
\[
H^k(X;\Z )=
\begin{cases}
\Z & k=0, 1, 4 \\
\Z /2\Z & k=2 \\
\Z \oplus \Z /2\Z  & k=3 \\
0  & \text{otherwise} . 
\end{cases}
\]
\end{ex}
\begin{ex}[Example \ref{ex3.8} with $k=0$]\label{ex6.6}
Let us calculate the cohomology for Example \ref{ex3.8} with $k=0$. In this case, $B$ 
is a compact, connected, oriented surface of genus $g\ge 0$, with boundary. (We 
include the case of $g=0$, in which case, $X=S^3\times S^1$.) 
Let $Q_{4g+1}$ be the polygon with $4g+1$ edges. $Q_{4g+1}$ 
gives a cell decomposition of $B$ with one zero-cell $e^{(0)}$, $2g+1$ one-cells 
$e^{(1)}_1$, $\ldots$, $e^{(1)}_{2g+1}$, and one two-cell $e^{(2)}$ 
as in Figure \ref{fig8}. 
\begin{figure}[htbp]
\begin{center}
\input{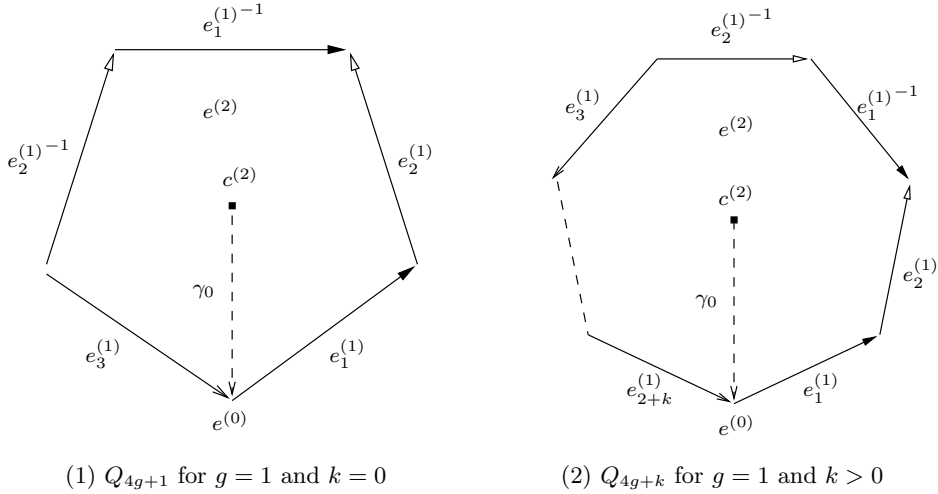}
\caption{polygons $Q_{4g+1}$ for $k=0$, $Q_{4g+k}$ for $k>0$, and the cell 
decompositions}
\label{fig8}
\end{center}
\end{figure}
In the cell decomposition, one-cell $e^{(1)}_{2i-1}$, $e^{(1)}_{2i}$, and $e^{(1)}_{2g+1}$ 
correspond to $\alpha_i$, $\beta_i$, and $\gamma$ in Figure \ref{fig2}, respectively. $B$ can be obtained 
from $Q_{4g+1}$ by identifying the oriented arrows corresponding to $e^{(1)}_i$ for $i=1$, 
$\ldots$, $2g$. Then we have the natural map $\varphi :Q_{4g+1}\to B$. 
In this case, the pull-back bundle $\varphi^*T^2_P$ is identified with the trivial bundle 
$Q_{4g+1}\times T^2$ and the pull-back $\varphi^*X$ of $X$ is identified with 
the quotient space of $Q_{4g+1}\times T^2$ which is obtained by collapsing each fiber 
on the one-cell $e^{(1)}_{2g+1}$ and all vertices of $Q_{4g+1}$ with $0\times S^1$. 
In this case, all fibers $\mu^{-1}(c^{(p)}_{\lambda})$ except for 
$\mu^{-1}(c^{(0)})$ and $\mu^{-1}(c^{(1)}_{2g+1})$ are diffeomorphic to 
$\pi_T^{-1}(c^{(p)}_{\lambda})\cong T^2$, whereas on the cell 
$e^{(p)}_{\lambda}=e^{(0)}$, $e^{(1)}_{2g+1}$, the fiber $\mu^{-1}(c^{(p)}_{\lambda})$ 
is diffeomorphic to $T^2/(0\times S^1)$ and the map 
$\nu_{c^{(p)}_{\lambda}}:\pi_T^{-1}(c^{(p)}_{\lambda})\to \mu^{-1}(c^{(p)}_{\lambda})$ 
can be identified with the natural projection $T^2\to T^2/(0\times S^1)$. 
For $e^{(0)}$ and $e^{(1)}_{2g+1}$, 
$\nu_{c^{(p)}_{\lambda}}^*\left( H^q(\mu^{-1}(c^{(p)}_{\lambda}); \Z )\right)$ 
can be obtained by (\ref{eq6.1}) in Example \ref{ex6.4}, 
and for the other cells, $\nu_{c^{(p)}_{\lambda}}^*:
H^q(\mu^{-1}(c^{(p)}_{\lambda}); \Z )\to 
H^q(\pi_T^{-1}(c^{(p)}_{\lambda}); \Z )$ is an isomorphism.  

For $q=0$, the cohomology $H^p(B ; {\mathcal H}_X^0)$ of 
$\left( C^p(B ; {\mathcal H}_X^0), \delta \right)$ is naturally identified with 
the cohomology $H^p(B ; \Z )$ with $\Z$-coefficient. 

For $q=1$, the degree $p$ cochain $u\in C^p(B ;{\mathcal H}^1_X)$ takes 
values as follows 
\[
u(e^{(p)}_{\lambda})\in \nu_{c^{(p)}_{\lambda}}^*\left( H^q(\mu^{-1}(c^{(p)}_{\lambda}); \Z )\right) =
\begin{cases}
\Z \oplus 0 & p=0\ \text{or}\ p=1\ \text{and}\ \lambda =2g+1\\ 
\Z \oplus \Z & p=1\ \text{and}\ \lambda =1, \ldots , 2g\ \text{or}\ p=2 \\
0 & \text{otherwise. }
\end{cases}
\] 
All differentials 
$\delta^p:C^p(B ;{\mathcal H}^1_X)\to C^{p+1}(B ;{\mathcal H}^1_X)$ vanish 
except for $p=0,1$, and in this case of $p=0$, $\delta^0$ is given as follows 
\begin{align}
\left( \delta^0u\right) (e^{(1)}_{2i-1})&= ^t\rho (\alpha_i)^{-1}u(e^{(0)})- u(e^{(0)})
\ \text{for}\ i=1, \ldots ,g, \nonumber \\
\left( \delta^0u\right) (e^{(1)}_{2i})&= ^t\rho (\beta_i)^{-1}u(e^{(0)})- u(e^{(0)}) 
\ \text{for}\ i=1, \ldots ,g, \nonumber \\
\left( \delta^0u\right) (e^{(1)}_{2g+1})&= ^t\rho (\gamma)^{-1}u(e^{(0)})- u(e^{(0)})
=0. 
\nonumber 
\end{align}
We identify each fiber $\pi_T^{-1}(c^{(1)}_{\lambda})$ with 
$\pi_T^{-1}(c^{(0)})$ as in Example \ref{ex6.4}. Then $\delta^1$ is given by  
\begin{align}
&\left( \delta^1u\right) (e^{(2)}) \nonumber \\
&= u(e^{(1)}_1)
				 + ^t\rho (\alpha_1)^{-1}u(e^{(1)}_2)
- ^t\rho (\alpha_1\beta_1\alpha_1^{-1})^{-1}u(e^{(1)}_1) 
- ^t\rho ([\alpha_1 \beta_1])^{-1}u(e^{(1)}_2)  \nonumber \\
&+\cdots \nonumber \\				  
&+ ^t\rho (\prod_{i=1}^{g-1}[\alpha_i \beta_i])^{-1}u(e^{(1)}_{2g-1})
+ ^t\rho (\prod_{i=1}^{g-1}[\alpha_i \beta_i]\alpha_g)^{-1}u(e^{(1)}_{2g})
\nonumber \\
&- ^t\rho (\prod_{i=1}^{g-1}[\alpha_i \beta_i]\alpha_g\beta_g\alpha_g^{-1})^{-1}
u(e^{(1)}_{2g-1})- ^t\rho (\prod_{i=1}^g[\alpha_i \beta_i])^{-1}u(e^{(1)}_{2g}) \nonumber \\
&+ ^t\rho (\prod_{i=1}^g[\alpha_i \beta_i])^{-1}u(e^{(1)}_{2g+1}) \nonumber \\
&=\left( 1- ^t\rho (\beta_1)^{-1}\right) u(e^{(1)}_1)+
\left( ^t\rho (\alpha_1)^{-1}-1\right) u(e^{(1)}_2)\nonumber \\
&+\cdots \nonumber \\
&=\left( 1- ^t\rho (\beta_g)^{-1}\right) u(e^{(1)}_{2g-1})+
\left( ^t\rho (\alpha_g)^{-1}-1\right) u(e^{(1)}_{2g})+u(e^{(1)}_{2g+1})  .
\nonumber 
\end{align}
Then the cohomology groups are calculated by 
\[
H^p(B ; {\mathcal H}_X^1)=
\begin{cases}
\Z & p=0, 2 \\
\Z^{\oplus 4g} & p=1 \\
0 & \text{otherwise}
\end{cases}
\]
for all $a_i=b_i=0$. In the other case, 
\[
H^p(B ; {\mathcal H}_X^1)=
\begin{cases}
\Lambda & p=1 \\
\Z/(a_i, b_j)\Z & p=2 \\
0 & \text{otherwise,}
\end{cases}
\]
where $(a_i, b_j)$ is the greatest common measure of all 
$a_i$ and $b_j$ which are not equal to $0$ and 
\[
\begin{split}
\Lambda =&\left\{ (u_1, v_1, \ldots , u_g, v_g)\in \Z^{\oplus 2g}\colon 
\sum_{i=1}^g(b_iu_i-a_iv_i)=0 \right\} \\
&\oplus \Z^{\oplus 2g}/
\left\{ (-a_1u, -b_1u, \ldots ,-a_gu, -b_gu)\colon u\in \Z \right\} . 
\end{split}
\] 

For $q=2$, the degree $p$ cochain $u\in C^p(B ;{\mathcal H}^2_X)$ takes 
values as follows 
\[
u(e^{(p)}_{\lambda})\in \nu_{c^{(p)}_{\lambda}}^*
\left( H^q(\mu^{-1}(c^{(p)}_{\lambda});\Z )\right) =
\begin{cases}
\Z  & p=1\ \text{and}\ \lambda =1, \ldots , 2g\ \text{or}\ p=2 \\
0 & \text{otherwise. }
\end{cases}
\]
All differentials 
$\delta^p:C^p(B ;{\mathcal H}^2_X)\to C^{p+1}(B ;{\mathcal H}^2_X)$ vanish. 
It is clear except for $p=1$. 
In the case of $p=1$, since all holonomies along $e^{(1)}_\lambda$ induce the identity of 
$H^2(\pi_T^{-1}(c^{(0)});\Z )$, the differential $\delta^1$ is given by  
\[
\begin{split}
\delta^1u(e^{(2)})=&u(e^{(1)}_1)+u(e^{(1)}_2)-u(e^{(1)}_1)-u(e^{(1)}_2) \\ 
&+\cdots \\
&+u(e^{(1)}_{2g-1})+u(e^{(1)}_{2g})-u(e^{(1)}_{2g-1})-u(e^{(1)}_{2g})+u(e^{(1)}_{2g+1})=0. 
\end{split}
\]
Then the cohomology groups are obtained by 
\[
H^p(B ; {\mathcal H}_X^2)=
\begin{cases}
\Z^{\oplus 2g} & p=1 \\
\Z & p=2 \\
0 & \text{otherwise. }
\end{cases}
\]
The table for the $E_2$-terms is in Figure \ref{fig9}. 
\begin{figure}[hbtp]
\begin{center}
\input{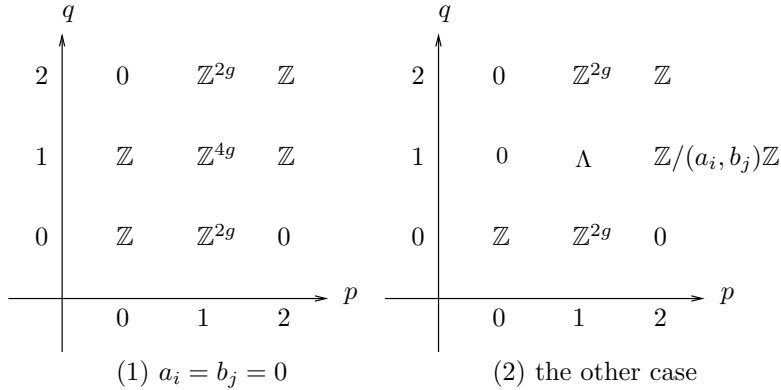}
\caption{the table of $(E_X^{p,q})_2$-terms for $k=0$ in Example \ref{ex3.8}}
\label{fig9}
\end{center}
\end{figure}
In particular, the Leray spectral sequence is 
degenerate at $E^2$-term, too, and the cohomology groups of $X$ are given by
\[
H^k(X;\Z )=
\begin{cases}
\Z & k=0, 4 \\
\Z^{\oplus 2g+1} & k=1, 3 \\
\Z^{\oplus 4g} & k=2 \\
0  & \text{otherwise} 
\end{cases}
\]
for all $a_i=b_i=0$. In the other case, 
\[
H^k(X;\Z )=
\begin{cases}
\Z & k=0, 4 \\
\Z^{\oplus 2g} & k=1 \\
\Lambda & k=2 \\
\Z^{\oplus 2g} \oplus \Z/(a_i, b_j)\Z & k=3 \\
0  & \text{otherwise.} 
\end{cases}
\]
\end{ex}
\begin{ex}[Example \ref{ex3.8} with $k=2$]
In the case of $k=2$ in Example \ref{ex3.8}, we must replace $Q_{4g+1}$ in Example 
\ref{ex6.6} with the polygon $Q_{4g+2}$ with $4g+2$ edges $e^{(1)}_1$, $e^{(1)}_2$, 
$(e^{(1)}_1)^{-1}$, $(e^{(1)}_2)^{-1}$, $\ldots$, $e^{(1)}_{2g-1}$, $e^{(1)}_{2g}$, 
$(e^{(1)}_{2g-1})^{-1}$, $(e^{(1)}_{2g})^{-1}$, $e^{(1)}_{2g+1}$, and 
$e^{(1)}_{2g+2}$. $e^{(1)}_{2i-1}$ and $e^{(1)}_{2i}$ correspond to $\alpha_i$, 
$\beta_i$ in Figure \ref{fig2}, respectively, and $e^{(1)}_{2g+1}$ and $e^{(1)}_{2g+2}$ correspond to 
the edge arcs $\gamma_1$ and $\gamma_2$ in Figure \ref{fig2}, respectively. 
See also Figure \ref{fig8}. 
As before, this gives $B$ a cell decomposition with two zero-cells $e^{(0)}_1$ and 
$e^{(0)}_2$, $2g+2$ one-cells $e^{(1)}_1$, $\ldots$, $e^{(1)}_{2g+2}$, and 
one two-cell $e^{(2)}$. By the same way in Example \ref{ex6.6}, we have the natural map 
$\varphi :Q_{4g+2}\to B$. In this case, the pull-back bundle $\varphi^*T^2_P$ 
is also identified with the trivial bundle $Q_{4g+2}\times T^2$ but the 
pull-back $\varphi^*X$ of $X$ is identified 
with its quotient space which is obtained from $Q_{4g+2}\times T^2$ by collapsing 
each fiber on the one-cell $e^{(1)}_{2g+1}$ with $S^1\times 0$, on the one-cell 
$e^{(1)}_{2g+2}$ with $0\times S^1$, and on all vertices of $Q_{4g+2}$ with $T^2$. 
The fibers $\mu^{-1}(c^{(p)}_{\lambda})$ are diffeomorphic as follows 
\[
\mu^{-1}(c^{(p)}_{\lambda})\cong 
\begin{cases}
\text{one point} & p=0 \\
T^2/S^1\times 0 & p=1\ \text{and}\ \lambda =2g+1 \\
T^2/0\times S^1 & p=1\ \text{and}\ \lambda =2g+2 \\
\pi_T^{-1}(c^{(p)}_{\lambda})\cong T^2 & \text{otherwise. }
\end{cases}
\]
For $e^{(0)}_{\lambda}$, $e^{(1)}_{2g+1}$, and $e^{(1)}_{2g+2}$,  
$\nu_{c^{(p)}_{\lambda}}^*\left( H^p(\mu^{-1}(c^{(p)}_{\lambda}); \Z )\right)$ 
is identified as follows
\[
\begin{split}
\nu_{c^{(0)}_{\lambda}}^*\left( H^q(\mu^{-1}(c^{(0)}_{\lambda}); \Z )\right)
&=
\begin{cases}
\Z & q=0 \\
0 & \text{otherwise, }
\end{cases} \\
\nu_{c^{(1)}_{2g+1}}^*\left( H^q(\mu^{-1}(c^{(1)}_{2g+1}); \Z )\right) &=
\begin{cases}
\Z & q=0 \\
0 \oplus \Z & q=1 \\
0 & \text{otherwise, }
\end{cases} \\
\nu_{c^{(1)}_{2g+2}}^*\left( H^q(\mu^{-1}(c^{(1)}_{2g+2}); \Z )\right) &=
\begin{cases}
\Z & q=0 \\
\Z \oplus 0 & q=1 \\
0 & \text{otherwise, }
\end{cases}
\end{split}
\]
and $\nu_{c^{(p)}_{\lambda}}^*:H^q(\mu^{-1}(c^{(p)}_{\lambda}); \Z )\to 
H^q(\pi_T^{-1}(c^{(p)}_{\lambda}); \Z )$ is isomorphic for the other cells.  
The similar calculus as in Example \ref{ex6.6} gives the table of the 
$E_2$-terms in Figure \ref{fig10}.  
The Leray spectral sequence is degenerate at $E^2$-term, and the cohomology 
groups of $X$ are given by
\[
H^k(X;\Z )=
\begin{cases}
\Z & k=0, 4 \\
\Z^{\oplus 2g} & k=1, 3 \\
\Z^{\oplus 4g} & k=2 \\
0  & \text{otherwise. } 
\end{cases}
\]
\end{ex}
\begin{ex}[Example \ref{ex3.8} with $k=3$]
In the case of $k=3$ in Example \ref{ex3.8}, we take $X_{\Delta}=\C P^2$ as 
a symplectic toric manifold with the triangle as its Delzant polytope. In this case, we change   
$Q_{4g+1}$ in Example \ref{ex6.6} with the polygon $Q_{4g+3}$ with $4g+3$ edges 
$e^{(1)}_1$, $e^{(1)}_2$, $(e^{(1)}_1)^{-1}$, $(e^{(1)}_2)^{-1}$, $\ldots$, 
$e^{(1)}_{2g-1}$, $e^{(1)}_{2g}$, $(e^{(1)}_{2g-1})^{-1}$, $(e^{(1)}_{2g})^{-1}$, 
$e^{(1)}_{2g+1}$, $\ldots$, $e^{(1)}_{2g+3}$. $e^{(1)}_{2i-1}$ and $e^{(1)}_{2i}$ 
correspond to $\alpha_i$, $\beta_i$ in Figure \ref{fig2}, respectively, and $e^{(1)}_{2g+1}$, $\ldots$, 
$e^{(1)}_{2g+3}$ correspond to the edge arcs $\gamma_1$, $\ldots$, 
$\gamma_3$ in Figure \ref{fig2}, respectively. See also Figure \ref{fig8}. This gives $B$ a cell 
decomposition with three zero-cells $e^{(0)}_1$, $\ldots$, $e^{(0)}_3$, $2g+3$ 
one-cells $e^{(1)}_1$, $\ldots$, $e^{(1)}_{2g+3}$, and one two-cell $e^{(2)}$. 
In this case, the pull-back bundle $\varphi^*T^2_P$ 
is also identified with the trivial bundle $Q_{4g+3}\times T^2$ but the 
pull-back $\varphi^*X$ of $X$ is identified 
with its quotient space which is obtained from $Q_{4g+3}\times T^2$ by collapsing 
each fiber on the one-cell $e^{(1)}_{2g+1}$ with $S^1\times 0$, on the one-cell 
$e^{(1)}_{2g+2}$ with $0\times S^1$, on the one-cell $e^{(1)}_{2g+3}$ with the circle 
in $T^2$ generated by $(-1, -1)$ which we denote by $-{\rm diag}(S^1)$ , and on all vertices 
of $Q_{4g+3}$ with $T^2$. The fibers $\mu^{-1}(c^{(p)}_{\lambda})$ are diffeomorphic as follows 
\[
\mu^{-1}(c^{(p)}_{\lambda})\cong 
\begin{cases}
\text{one point} & p=0 \\
T^2/S^1\times 0 & p=1\ \text{and}\ \lambda =2g+1 \\
T^2/0\times S^1 & p=1\ \text{and}\ \lambda =2g+2 \\
T^2/-{\rm diag}(S^1) & P=1\ \text{and}\ \lambda =2g+3 \\
\pi_T^{-1}(c^{(p)}_{\lambda})\cong T^2 & \text{otherwise. }
\end{cases}
\]
For $e^{(0)}_{\lambda}$, $e^{(1)}_{2g+1}$, $\ldots$, $e^{(1)}_{2g+3}$,  
$\nu_{c^{(p)}_{\lambda}}^*\left( H^p(\mu^{-1}(c^{(p)}_{\lambda}); \Z )\right)$ 
is identified as follows
\[
\begin{split}
\nu_{c^{(0)}_{\lambda}}^*\left( H^q(\mu^{-1}(c^{(0)}_{\lambda}); \Z )\right) &=
\begin{cases}
\Z & q=0 \\
0 & \text{otherwise, }
\end{cases} \\
\nu_{c^{(1)}_{2g+1}}^*\left( H^q(\mu^{-1}(c^{(1)}_{2g+1}); \Z )\right) &=
\begin{cases}
\Z & q=0 \\
0 \oplus \Z & q=1 \\
0 & \text{otherwise, }
\end{cases} \\
\nu_{c^{(1)}_{2g+2}}^*\left( H^q(\mu^{-1}(c^{(1)}_{2g+2}); \Z )\right) &=
\begin{cases}
\Z & q=0 \\
\Z \oplus 0 & q=1 \\
0 & \text{otherwise, }
\end{cases} \\
\nu_{c^{(1)}_{2g+3}}^*\left( H^q(\mu^{-1}(c^{(1)}_{2g+3}); \Z )\right) &=
\begin{cases}
\Z & q=0 \\
{\rm offdiag}(\Z ) & q=1 \\
0 & \text{otherwise, }
\end{cases}
\end{split}
\]
where ${\rm offdiag}(\Z )$ is the sub-lattice of $\Z \oplus \Z$ which is 
generated by $(1, -1)$. For other cells, 
$\nu_{c^{(p)}_{\lambda}}^*: H^p(\mu^{-1}(c^{(p)}_{\lambda}); \Z )\to
H^p(\pi_T^{-1}(c^{(p)}_{\lambda}); \Z )$ is isomorphic.  
The similar calculus as in Example \ref{ex6.6} gives the table of the 
$E_2$-terms in Figure \ref{fig10}.  
Then the Leray spectral sequence is degenerate at $E^2$-term, and the cohomology groups 
of $X$ are given by
\[
H^k(X;\Z )=
\begin{cases}
\Z & k=0, 4 \\
\Z^{\oplus 2g} & k=1, 3 \\
\Z^{\oplus 4g+1} & k=2 \\
0  & \text{otherwise. } 
\end{cases}
\]
\end{ex}
\begin{ex}[Example \ref{ex3.8} with $k=4$]
In the case of $k=4$ in Example \ref{ex3.8}, we take 
$X_{\Delta}=S^2\times S^2$ as a symplectic toric manifold with the square 
as its Delzant polytope. In this case, we change   
$Q_{4g+1}$ in Example \ref{ex6.6} with the polygon $Q_{4g+4}$ with $4g+4$ edges 
$e^{(1)}_1$, $e^{(1)}_2$, $(e^{(1)}_1)^{-1}$, $(e^{(1)}_2)^{-1}$, $\ldots$, 
$e^{(1)}_{2g-1}$, $e^{(1)}_{2g}$, $(e^{(1)}_{2g-1})^{-1}$, $(e^{(1)}_{2g})^{-1}$, 
$e^{(1)}_{2g+1}$, $\ldots$, $e^{(1)}_{2g+4}$. $e^{(1)}_{2i-1}$ and $e^{(1)}_{2i}$ 
correspond to $\alpha_i$, $\beta_i$ in Figure \ref{fig2}, respectively, and $e^{(1)}_{2g+1}$, $\ldots$, 
$e^{(1)}_{2g+4}$ correspond to the edge arcs $\gamma_1$, $\ldots$, 
$\gamma_4$ in Figure \ref{fig2}, respectively. See also Figure \ref{fig8}. This gives $B$ a cell 
decomposition with four zero-cells $e^{(0)}_1$, $\ldots$, $e^{(0)}_4$, $2g+4$ 
one-cells $e^{(1)}_1$, $\ldots$, $e^{(1)}_{2g+4}$, and one two-cell $e^{(2)}$. 
In this case, the pull-back bundle $\varphi^*T^2_P$ 
is also identified with the trivial bundle $Q_{4g+4}\times T^2$ but the 
pull-back $\varphi^*X$ of $X$ is identified with its quotient space which is 
obtained from $Q_{4g+4}\times T^2$ by collapsing 
each fiber on the one-cell $e^{(1)}_{2g+1}$ with $S^1\times 0$, on the one-cell 
$e^{(1)}_{2g+2}$ with $0\times S^1$, on the one-cell $e^{(1)}_{2g+3}$ with 
$-S^1\times 0$, on the one-cell $e^{(1)}_{2g+4}$ with $0\times -S^1$, and on all 
vertices of $Q_{4g+4}$ with $T^2$, where $-S^1\times 0$ is the circle generated by 
$(-1,0)$ etc. The fibers $\mu^{-1}(c^{(p)}_{\lambda})$ are diffeomorphic as follows 
\[
\mu^{-1}(c^{(p)}_{\lambda})\cong 
\begin{cases}
\text{one point} & p=0 \\
T^2/S^1\times 0 & p=1\ \text{and}\ \lambda =2g+1 \\
T^2/0\times S^1 & p=1\ \text{and}\ \lambda =2g+2 \\
T^2/-S^1\times 0 & P=1\ \text{and}\ \lambda =2g+3 \\
T^2/0\times -S^1 & P=1\ \text{and}\ \lambda =2g+4 \\
\pi_T^{-1}(c^{(p)}_{\lambda})\cong T^2 & \text{otherwise. }
\end{cases}
\]
For $e^{(0)}_{\lambda}$, $e^{(1)}_{2g+1}$, $\ldots$, $e^{(1)}_{2g+4}$,  
$\nu_{c^{(p)}_{\lambda}}^*\left( H^p(\mu^{-1}(c^{(p)}_{\lambda}); \Z )\right)$ 
is identified as follows
\[
\begin{split}
\nu_{c^{(0)}_{\lambda}}^*\left( H^q(\mu^{-1}(c^{(0)}_{\lambda}); \Z )\right) &=
\begin{cases}
\Z & q=0 \\
0 & \text{otherwise, }
\end{cases} \\
\nu_{c^{(1)}_{\lambda}}^*\left( H^q(\mu^{-1}(c^{(1)}_{\lambda}); \Z )\right) &=
\begin{cases}
\Z & q=0 \\
0 \oplus \Z & q=1 \\
0 & \text{otherwise, }
\end{cases}\ (\lambda =2g+1, 2g+3)\\
\nu_{c^{(1)}_{\lambda}}^*\left( H^q(\mu^{-1}(c^{(1)}_{\lambda}); \Z )\right) &=
\begin{cases}
\Z & q=0 \\
\Z \oplus 0 & q=1 \\
0 & \text{otherwise}
\end{cases}\ (\lambda =2g+2, 2g+4), 
\end{split}
\]
and $\nu_{c^{(p)}_{\lambda}}^*: H^p(\mu^{-1}(c^{(p)}_{\lambda}); \Z )\to 
H^p(\pi_T^{-1}(c^{(p)}_{\lambda}); \Z )$ is an isomorphism for the other 
cells. The similar calculus as in Example \ref{ex6.6} gives the table of the 
$E_2$-terms in Figure \ref{fig10}.  
\begin{figure}[hbtp]
\begin{center}
\input{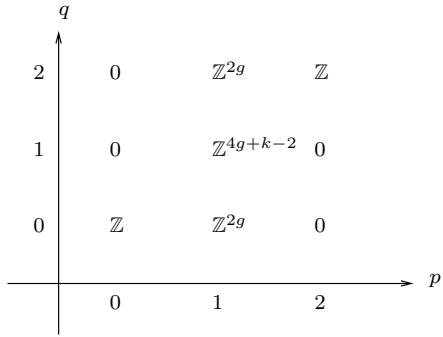}
\caption{$(E_X^{p,q})_2$-terms for $k=2,3,4$ in Example \ref{ex3.8}}
\label{fig10}
\end{center}
\end{figure}
Then the Leray spectral sequence is degenerate at $E^2$-term, and the cohomology groups 
of $X$ are given by
\[
H^k(X;\Z )=
\begin{cases}
\Z & k=0, 4 \\
\Z^{\oplus 2g} & k=1, 3 \\
\Z^{\oplus 4g+2} & k=2 \\
0  & \text{otherwise. } 
\end{cases}
\]
\end{ex}
\begin{ex}[Example \ref{ex3.9}]
Let us calculate the cohomology groups for Example \ref{ex3.9}. In this case, we take    
the pentagon $Q_5$ with edges $e^{(1)}_1$, $e^{(1)}_2$, $(e^{(1)}_1)^{-1}$, 
$(e^{(1)}_2)^{-1}$, and $e^{(1)}_3$. $e^{(1)}_1$, $e^{(1)}_2$, and $e^{(1)}_3$ correspond to 
$\alpha$, $\beta$, and the edge arc $\gamma$ in Figure \ref{fig2}, respectively. 
See also Figure \ref{fig8}. This gives $B$ a cell decomposition with one zero-cell $e^{(0)}$, 
three one-cells $e^{(1)}_1$, $\ldots$, $e^{(1)}_3$, and one two-cell $e^{(2)}$. 
In this case, the pull-back bundle $\varphi^*T^2_P$ is also identified with the trivial 
bundle $Q_5\times T^2$ but the pull-back $\varphi^*X$ of $X$ is identified with its quotient space which 
is obtained from $Q_5\times T^2$ by collapsing each fiber on the one-cell 
$e^{(1)}_3$ with $0\times S^1$, and on all vertices of $Q_5$ with $T^2$. 
The fibers $\mu^{-1}(c^{(p)}_{\lambda})$ are diffeomorphic as follows 
\[
\mu^{-1}(c^{(p)}_{\lambda})\cong 
\begin{cases}
\text{one point} & p=0 \\
T^2/0\times S^1 & p=1\ \text{and}\ \lambda =3 \\
\pi_T^{-1}(c^{(p)}_{\lambda})\cong T^2 & \text{otherwise. }
\end{cases}
\]
For $e^{(0)}$ and $e^{(1)}_3$, 
$\nu_{c^{(p)}_{\lambda}}^*\left( H^p(\mu^{-1}(c^{(p)}_{\lambda}); \Z )\right)$ 
is identified as follows
\[
\begin{split}
\nu_{c^{(0)}}^*\left( H^q(\mu^{-1}(c^{(0)}); \Z )\right) &=
\begin{cases}
\Z & q=0 \\
0 & \text{otherwise, }
\end{cases} \\
\nu_{c^{(1)}_3}^*\left( H^q(\mu^{-1}(c^{(1)}_3); \Z )\right) &=
\begin{cases}
\Z & q=0 \\
\Z \oplus 0 & q=1 \\
0 & \text{otherwise. }
\end{cases}
\end{split}
\]
and $\nu_{c^{(p)}_{\lambda}}^*: H^p(\mu^{-1}(c^{(p)}_{\lambda}); \Z )\to 
H^p(\pi_T^{-1}(c^{(p)}_{\lambda}); \Z )$ is an isomorphism for the other 
cells. The similar calculus as in Example \ref{ex6.6} gives the table of 
the $E_2$-terms in Figure \ref{fig11}.  
\begin{figure}[hbtp]
\begin{center}
\input{fig10.pstex_t}
\caption{$(E_X^{p,q})_2$-terms for Example \ref{ex3.9}}
\label{fig11}
\end{center}
\end{figure}
Then the Leray spectral sequence is degenerate at $E^2$-term, and the cohomology groups 
of $X$ are given by
\[
H^k(X;\Z )=
\begin{cases}
\Z & k=0, 4 \\
\Z^{\oplus 2} & k=1, 3 \\
\Z^{\oplus 3} & k=2 \\
0  & \text{otherwise. } 
\end{cases}
\]
\end{ex}

\subsection{Signatures}
In this subsection, we shall give the method of computing the signature for a 
four-dimensional case by using 
the Novikov additivity. Let $B$ be a surface with at least one corner, and 
$X$ a twisted toric manifold associated with a principal $SL_2(\Z )$-bundle 
$P$ on $B$. For simplicity, assume that $B$ has only one boundary component. 
We divide $B$ into two parts $B_1$ and $B_2$, where $B_2$ is the closed 
neighborhood of the boundary $\partial B$ such that $\partial B$ is a deformation 
retract of $B_2$ and $B_1$ is the closure $B_1=\overline{B\backslash B_2}$ of 
the remainder. 
We set $X_i=\mu^{-1}(B_i)$ for $i=1,2$, and denote by $\sigma (X_i)$ and 
$\sigma (X)$ the signature of $X_i$ and $X$, respectively. The Novikov additivity says that
\begin{equation}\label{signX}
\sigma (X)=\sigma (X_1)+\sigma (X_2). 
\end{equation}

First let us compute the signature $\sigma (X_1)$ of $X_1$. We notice that $X_1$ 
is the associated $T^2$-bundle for $P$. When the genus of $B$ is equal to zero, $B_1$ 
is contractible. In this case, 
the signature $\sigma (X_1)$ is zero. 

\begin{figure}[hbtp]
\begin{center}
\input{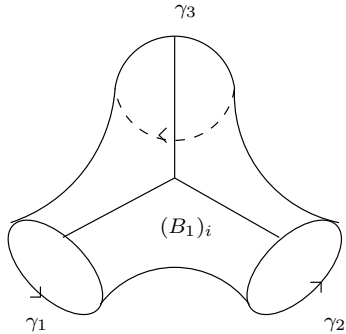}
\caption{$(B_1)_i$ and $\gamma_j$}
\label{fig13}
\end{center}
\end{figure}

When the genus of $B$ is greater than zero, we give $B_1$ a trinion decomposition 
$\displaystyle B_1=\cup_{i=1}^k(B_1)_i$, where each $(B_1)_i$ is a trinion, that is, 
a surface obtained from $S^2$ by removing three distinct open discs. Let 
$(X_1)_i=\mu^{-1}((B_1)_i)$ for $i=1, \ldots ,k$. From the Novikov additivity, we have 
\begin{equation}\label{signX_1}
\sigma (X_1)=\sum_{i=1}^k\sigma ((X_1)_i) ,  
\end{equation}
and each $\sigma ((X_1)_i)$ can be computed as follows. 
We take the oriented boundary loops $\gamma_1$, $\gamma_2$, and $\gamma_3$ of $(B_1)_i$ 
as in Figure \ref{fig13} which represent generators of $\pi_1((B_1)_i)$ with 
$[\gamma_1]\cdot [\gamma_2]\cdot [\gamma_3]=1$.   
Let $\rho :\pi_1((B_1)_i)\to SL_2(\Z )(=Sp(2;\Z ))$ be the representation which 
determines $T^2$-bundle $(X_1)_i$ on $(B_1)_i$. We set $C_j=\rho ([\gamma_j])$ for 
$j=1, 2, 3$. 
For $C_1$ and $C_2$, define the vector space $V_{C_1, C_2}$ and the bilinear form 
$\< \ ,\ \>_{C_1, C_2}$ on $V_{C_1, C_2}$ by 
\[
\begin{split}
&V_{C_1, C_2}=\{ (x, y)\in \R^2\times \R^2 \colon (C_1^{-1}-I)x+(C_2-I)y=0\} , \\
&\< (x, y), (x', y')\>_{C_1, C_2}= ^t(x+y)J(I-C_2)y'
\end{split}
\]
for $(x, y)$, $(x', y')\in V_{C_1, C_2}$, where $I=
\begin{pmatrix}
1 & 0 \\
0 & 1
\end{pmatrix}$ and 
$J=
\begin{pmatrix}
0 & 1 \\
-1 & 0
\end{pmatrix}$. 
It is easy to show that $\< \ ,\ \>_{C_1, C_2}$ is symmetric and we denote the signature 
of $\< \ ,\ \>_{C_1, C_2}$ by $\tau_1(C_1, C_2)$. 
\begin{thm}[\cite{E, Me}]\label{sign(X_1)_i}
$\sigma ((X_1)_i)=\tau_1(C_1, C_2)$. 
\end{thm}
We should notice that our orientation of $(X_1)_i$ is different from that in \cite{E, Me}. 
From (\ref{signX_1}) and Theorem \ref{sign(X_1)_i}, we can compute $\sigma (X_1)$.  
\begin{rem}
Meyer shows in \cite{Me} that $\tau_1$ defines the cocycle 
$\tau_1:Sp(2;\Z )\times Sp(2;\Z )\to \Z$ of $Sp(2;\Z )$ which is called Meyer's signature 
cocycle. The author was taught Meyer's signature cocycle by Endo \cite{E}. 
\end{rem}

Next, we shall compute $\sigma (X_2)$. For $X_2$, we can show the following lemma.  
\begin{lem}
$\mu^{-1}(\partial B)$ is a deformation retract of $X_2$.  
\end{lem}
\begin{proof}
Let $h:B_2\times I\to B_2$ be a deformation retraction with 
$h(\cdot , 0)=\id_{B_2}$ and $h(\cdot ,1)\in \partial B$. Define the map 
$\widetilde{h}:X_2\times I\to X_2$ of $h$ by 
\[
\widetilde{h}(x, s)=\nu \circ {\rm Hol}_{\gamma_{\mu(x), s}}(x)
\] 
for $(x, s)\in X_2\times I$, where $\gamma_{\mu(x), s}$ is the path 
$\gamma_{\mu (x),s}:I\to B_2$ which is defined by 
$\gamma_{\mu (x),s}(t)=h(\mu(x), st)$ for $t\in I$ and 
${\rm Hol}_{\gamma_{\mu(x), s}}$ is the parallel transport 
${\rm Hol}_{\gamma_{\mu(x), s}}:\pi_T^{-1}(\mu(x))\to \pi_T^{-1}
(\gamma_{\mu(x),s}(1))$ of $T^2_P$ along $\gamma_{\mu(x),s}$ with 
respect to the connection induced from that of $P$. Then $\widetilde{h}$ is 
a deformation retraction. 
\end{proof}

Assume that $B$ has $k$ corner points. Then the one-dimensional strata 
$\CS^{(1)}B$ has exactly $k$ connected components $(\CS^{(1)}B)_1$, 
$\ldots$, $(\CS^{(1)}B)_k$, and $\displaystyle \mu^{-1}(\partial B)
=\cup_{i=1}^k\mu^{-1}(\overline{(\CS^{(1)}B)_i})$, where 
$\overline{(\CS^{(1)}B)_i}$ is the closure of $(\CS^{(1)}B)_i$. 
It is easy to see from 
the similar construction of the twisted toric structure on a neighborhood of 
$\partial B$ as in Example \ref{ex3.9} and Theorem \ref{classification}  
that each $\mu^{-1}(\overline{(\CS^{(1)}B)_i})$ is homeomorphic to the 
two-dimensional sphere $S^2$ if $k\ge 2$, and is homeomorphic to $S^2$ 
with one self-intersection at north and south points if $k=1$.  
If $\overline{(\CS^{(1)}B)_i}\cap \overline{(\CS^{(1)}B)_j}\neq \emptyset$ 
for $i\neq j$, then they have two intersections if $k=2$, and have one intersection 
if $k>2$. In all cases, the intersections are transversal since by definition, 
a neighborhood of each intersection in $X$ is identified with that 
of the intersection of $\C\times \{ 0\}$ and $\{ 0\}\times \C$ in $\C^2$. 
Then $\mu^{-1}(\partial B)$ looks like a necklace of $k$ spheres and 
the homology group of $X_2$ is given by 
\[
H_p(X_2;\Z )=H_p(\mu^{-1}(\partial B); \Z )=
\begin{cases}
\Z & p=0,1 \\
\Z^{\oplus k} & p=2 \\
0 & \text{otherwise.}
\end{cases}
\]
Moreover the homology classes $[\mu^{-1}(\overline{(\CS^{(1)}B)_1})]$, $\ldots$, 
$[\mu^{-1}(\overline{(\CS^{(1)}B)_k})] \in H_2(X_2;\Z )$ 
represented by $\mu^{-1}(\overline{(\CS^{(1)}B)_i})$ for $i=1, \ldots , k$ 
are generators of $H_2(X_2;\Z )$. We set $S^2_i=\mu^{-1}(\overline{(\CS^{(1)}B)_i})$ 
for $i=1, \ldots ,k$. From the above fact, we can obtain the following proposition. 
\begin{prop}\label{int}
For $i\neq j$, the intersection number $[S^2_i]\cdot [S^2_j]$ of $[S^2_i]$ and 
$[S^2_j]$ is given as follows 
\[
[S^2_i]\cdot [S^2_j]=
\begin{cases}
0 & S^2_i\cap S^2_j=\emptyset \\
1 & S^2_i\cap S^2_j\neq \emptyset \ \text{and}\ k>2 \\
2 & S^2_i\cap S^2_j\neq \emptyset \ \text{and}\ k=2. 
\end{cases}
\]
\end{prop}
Assume that $k>1$. Let us compute the self-intersection number of $[S^2_i]$. 
We can take a contractible neighborhood $U$ of $\overline{(\CS^{(1)}B)_i}$ in $B$ so that 
$U\cap \CS^{(1)}B$ has exactly two connected components $(U\cap \CS^{(1)}B)_1$ 
and $(U\cap \CS^{(1)}B)_2$ except for $(\CS^{(1)}B)_i$. We may assume that 
$(U\cap \CS^{(1)}B)_1$ and $(U\cap \CS^{(1)}B)_2$ are located in Figure \ref{fig12}.    
\begin{figure}[hbtp]
\begin{center}
\input{sigma18.pstex_t}
\caption{$(\CS^{(1)}B)_i$ and $(U\cap \CS^{(1)}B)_a$}
\label{fig12}
\end{center}
\end{figure}
Let $\pi_{\CL}:\CL \to \CS^{(1)}B$ be the characteristic bundle of $X$. 
Since $\CL$ is primitive, we can take the local trivialization 
$\varphi^{\Z}:\pi_{\Z}^{-1}(U)\cong U\times \Z^2$ of $\pi_{\Z}:\Z^2_P\to B$ 
so that $\varphi^{\Z}$ also gives the trivializations 
$\pi_{\CL}^{-1}((\CS^{(1)}B)_i)\cong (\CS^{(1)}B)_i\times L$ and 
$\pi_{\CL}^{-1}((U\cap \CS^{(1)}B)_a)\cong (U\cap \CS^{(1)}B)_a\times L_a$ 
on $(\CS^{(1)}B)_i$ and $(U\cap \CS^{(1)}B)_a$ for $a=1, 2$, where $L$ and 
$L_a$ are rank one sub-lattices of $\Z^2$. We take the generators $u$, $u_a$ of 
$L$ and $L_a$ such that both of the determinants of $(u_1, u)$ and $(u, u_2)$ 
are equal to one, where $(u_1, u)$ (resp. $(u, u_2)$) denotes the matrix given by 
arranging the column vectors $u_1$ and $u$ (resp. $u$ and $u_2$) in this order. 
\begin{prop}\label{selfint}
The self-intersection number $[S^2_i]\cdot [S^2_i]$ is equal to the negative 
determinant $-\det (u_1, u_2)$ of $(u_1, u_2)$. 
\end{prop}
\begin{rem}
The determinant of $(u_1, u_2)$ does not depend on the choice of the local trivialization 
$\varphi^{\Z}$ since the structure group of the bundle is $SL_2(\Z )$.  
\end{rem}
\begin{proof}
Since the self-intersection number of $S^2_i$ is equal to the Euler number 
$\displaystyle \int_{S^2_i}e(\CN_{S^2_i})$ of the normal bundle $\CN_{S^2_i}$ 
of $S^2_i$ in $X$, for example, see \cite{BT}, we identify $\CN_{S^2_i}$. 
If necessary, by replacing the local trivialization, we may assume that $u_1=
\begin{pmatrix}
1\\
0
\end{pmatrix}$ 
and $u=
\begin{pmatrix}
0\\
1
\end{pmatrix}$. 
Then $u_2$ must be of the form $u_2= 
\begin{pmatrix}
-1\\
m
\end{pmatrix}$ 
for some $m\in \Z$ since $\det (u, u_2)=1$. The primitive vector $u_2$ 
defines the Hamiltonian $S^1$-action on $(\C^2, \om_{\C^2})$ by 
\[
t\cdot z=(e^{-2\pi \sqrt{-1}t}z_1, e^{2\pi \sqrt{-1}mt}z_2) 
\]
with the moment map 
\[
\mu_{u_2}(z)=-\lvert z_1\rvert^2+m\lvert z_2\rvert^2. 
\]
For $\varepsilon<0$, $\overline{\C^2}_{\mu_{u_2}\ge \varepsilon}$ denotes the 
cut space which is obtained from $\C^2$ by the symplectic cutting for 
this circle action. More precisely, let $\Phi :\C^2\times \C \to \R$ 
be the map defined by 
\[
\Phi (z, w)=-\abs{z_1}^2+m\abs{z_2}^2-\abs{w}^2-\varepsilon 
\]
for $(z, w)\in \C^2\times \C$. By the construction of the symplectic cutting, 
$\overline{\C^2}_{\mu_{u_2}\ge \varepsilon}=\Phi^{-1}(0)/S^1$, where the 
$S^1$-action on $\Phi^{-1}(0)$ is 
\[
t\cdot (z, w)=(t\cdot z, e^{-2\pi \sqrt{-1}t}w). 
\]
Then it is easy to see from the similar construction of the twisted toric 
structure on a neighborhood of $\partial B$ as in Example \ref{ex3.9} 
and Theorem \ref{classification} that the the neighborhood of $S^2_i$ 
in $X$ is orientation preserving diffeomorphic to that of 
\[
\{ (z, w)\in \Phi^{-1}(0)\colon z_2=0\} /S^1
\] 
in $\overline{\C^2}_{\mu_{u_2}\ge \varepsilon}$ for some $\varepsilon <0$. 
By this diffeomorphism, we identify $S^2_i$ with 
\[
\{ (z, w)\in \Phi^{-1}(0)\colon z_2=0\} /S^1. 
\]
Let $\pi_1 :\Phi^{-1}(0)\to \overline{\C^2}_{\mu_{u_2}\ge \varepsilon}$ and 
$\pi_2: \{ (z, w)\in \Phi^{-1}(0)\colon z_2=0\} \to S^2_i$ be the natural projections. 
It is easy to see from the symplectic cutting construction that 
\[
\begin{split}
&\pi_1^*T\overline{\C^2}_{\mu_{u_2}\ge \varepsilon}\oplus 
\underline{\Lie (S^1)\otimes_{\R}\C }
=T(\C^2\times \C)\Bigl|_{\Phi^{-1}(0)} \\
& \pi_2^*TS^2_i\oplus \underline{\Lie (S^1)\otimes_{\R}\C }
=T(\C\times \{0\} \times \C)\Bigl|_{\{(z, w)\in \Phi^{-1}(0)\colon z_2=0\} },   
\end{split}
\] 
where $\underline{\Lie (S^1)\otimes_{\R}\C}$ is the trivial complex line bundle. 
By the above identification, the normal bundle 
$\CN_{S^2_i}$ of $S^2_i$ in $X$ is isomorphic to the associated complex 
line bundle
\[  
\{(z, w)\in \Phi^{-1}(0)\colon z_2=0\} \times_{S^1} \C \to S^2_i 
\]
of $\pi_2: \{ (z, w)\in \Phi^{-1}(0)\colon z_2=0\} \to S^2_i$ with 
respect to the irreducible $S^1$-representation of weight $m$. 
This is the complex line bundle $\CO (-m)$ on $\C P^1$. Since the Euler 
class is equal to the first Chern class for a complex line bundle, 
the Euler number $\displaystyle \int_{S^2_i}e(\CN_{S^2_i})$ is equal 
to $-m$. 
\end{proof}
From Proposition \ref{int} and Proposition \ref{selfint}, we can compute the signature 
$\sigma (X_2)$ case-by-case for $k>1$. In the case of $k=1$, by blowing 
up the fiber on the corner point, which consists of the one point,
of $B$, we can reduce to the case of $k>1$ as in the following example. 
\begin{ex}
Let us compute the signature of the twisted toric manifold $X$ in Example 
\ref{ex3.9}. 
Recall that $B$ is a surface of genus one with one corner. As described above, 
we divide $B$ into two part $B_1$ and $B_2$. We give $B_1$ the trinion 
decomposition as in Figure \ref{fig14}. Then the easy computation shows 
that the value $\tau_1 (\rho ([\alpha^{-1}]), \rho ([\gamma^{-1}]))$ 
of the Meyer cocycle vanishes. This implies the signature $\sigma (X_1)$ 
of $X_1$ is zero. 
\begin{figure}[hbtp]
\begin{center}
\input{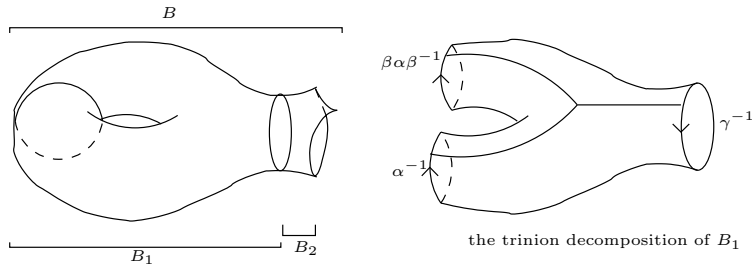}
\caption{$B$, $B_i$, and the trinion decomposition of $B_1$}
\label{fig14}
\end{center}
\end{figure}

Next we consider $X_2$. In general, the fiber of $\mu :X\to B$ on the corner point 
consists of only one point which we denote by $x_0$, and by definition, the 
neighborhood of $x_0$ is orientation preserving diffeomorphic to that of the origin of 
$\C^2$. Then we can blow up $X_2$ at $x_0$ and denote by 
$\widetilde{X_2}$ its blow-up. $\widetilde{X_2}$ is a twisted toric manifold on 
the surface $\widetilde{B_2}$ with two corners. For a blowing up of a symplectic 
toric manifold, see \cite{G, MS}. Let 
$(\CS^{(1)}\widetilde{B_2})_1$ and $(\CS^{(1)}\widetilde{B_2})_2$ denote 
connected components of the codimension one strata of $\widetilde{B_2}$ as in 
Figure \ref{fig15}.  
\begin{figure}[hbtp]
\begin{center}
\input{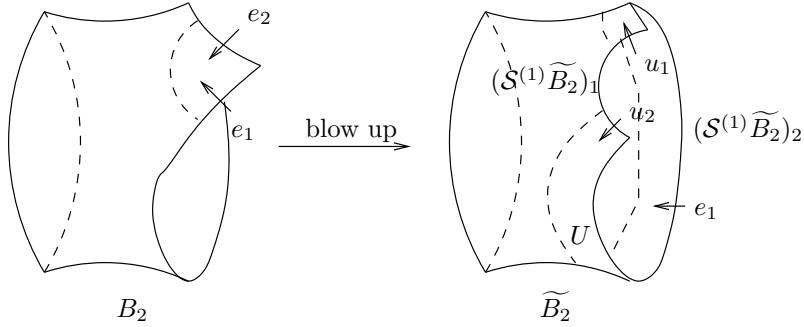}
\caption{$B_2$ and $\widetilde{B_2}$}
\label{fig15}
\end{center}
\end{figure}
Since the inverse image $S^2_1=\mu^{-1}(\overline{(\CS^{(1)}\widetilde{B_2})_1})$ of 
$(\CS^{(1)}\widetilde{B_2})_1$ is an exceptional divisor, its self-intersection 
number $[S^2_1]\cdot [S^2_1]$ is equal to $-1$. 

We compute the self-intersection number $[S^2_2]\cdot [S^2_2]$ of 
$S^2_2=\mu^{-1}(\overline{(\CS^{(1)}\widetilde{B_2})_2})$. By the construction of 
$X_2$ in Example \ref{ex3.9}, we can take $u_1$, $u_2$ in Proposition \ref{selfint} of 
the forms 
$u_1=
\begin{pmatrix}
4 \\
-1
\end{pmatrix}=
\begin{pmatrix}
3 & 1 \\
-1 & 0
\end{pmatrix}
\begin{pmatrix}
1 \\
1
\end{pmatrix}
$, 
$u_2=
\begin{pmatrix}
1 \\
1
\end{pmatrix}$, hence $[S^2_2]\cdot [S^2_2]=-\det (u_1, u_2)=-5$. 
The above computation and Proposition \ref{int} for $k=2$ show that the intersection 
matrix of $\widetilde{X_2}$ is 
$
\begin{pmatrix}
-1 & 2 \\
2 & -5
\end{pmatrix}$ and the signature $\sigma (\widetilde{X_2})$ is equal to $-2$. 
Since $\widetilde{X_2}$ is the blow-up of $X_2$ at $x_0$, the signature of $X_2$ is 
$\sigma (X_2)=\sigma (\widetilde{X_2})+1=-1$. Then the computations of 
$\sigma (X_1)$ and $\sigma (X_2)$ together with $(\ref{signX})$ shows 
that $\sigma (X)=-1$. 
\end{ex}

\appendix 
\section{Smoothness of cut spaces}
For $u_1, \ldots , u_d\in \Z^n$, define the Hamiltonian $T^d$-action 
on $(T^*T^n\times \C^d, \om_{T^*T^n}\oplus \om_{\C^d})$ by 
\begin{equation*}
t\cdot (\xi, \theta , z)=(\xi, \theta +\sum_{i=1}^dt_iu_i, (e^{-2\pi \i t_i}z_i)) . 
\end{equation*}
for $t=(t_1,\ldots ,t_d)\in T^d$ and $(\xi, \theta , z)\in T^*T^n\times \C^d$ 
with the moment map $\Phi :T^*T^n\times \C^d\to \R^d$ 
\[
\Phi (\xi, \theta , z)=(\< u_i,\xi \>-\abs{z_i}^2) .
\]
For $\lambda =(\lambda_1, \ldots , \lambda_d)\in \R^d$, we define the 
simultaneous cut space 
$\overline{(T^*T^n)}_{\{ \mu_{u_i}\ge \lambda_i \}_{i=1, \ldots ,d}}$ by 
\[ 
\overline{(T^*T^n)}_{\{ \mu_{u_i}\ge \lambda_i \}_{i=1, \ldots ,d}}
=\Phi^{-1}(\lambda )/T^d . 
\]
Let us investigate the smoothness of  
$\overline{(T^*T^n)}_{\{ \mu_{u_i}\ge \lambda_i \}_{i=1, \ldots ,d}}$. 
\begin{lem}\label{a1}
For a linear independent tuple $\{ u_1, \ldots , u_k\}$ of 
$k$ vectors in $\Z^n$, the following conditions are equivalent. 
\begin{enumerate}
\item $\{ u_1, \ldots , u_k\}$ is primitive
\item if $(t_1, \ldots , t_k)\in \R^k$ satisfies 
$\displaystyle \sum_{i=1}^kt_iu_i\in \Z^n$, then 
$(t_1, \ldots , t_k)\in \Z^k$
\end{enumerate}
\end{lem}
\begin{proof}
Let $A=(u_1, \ldots , u_k)\in M_{n\times k}(\Z )$, and 
$\varphi_A: \R^k\to \R^n$ be the linear map defined by $A$. 
Then 
\begin{align}
\text{(ii)}& \Longleftrightarrow \ \varphi_A^{-1}(\Z^n)=\Z^k \nonumber \\
                              & \Longleftrightarrow \ 
							  {\rm Im}\varphi_A\cap \Z^n=\varphi_A(\Z^k)\ \nonumber 
\end{align}
since $\varphi_A$ is injective. (The condition 
${\rm Im}\varphi_A\cap \Z^n\supset \varphi_A(\Z^k)$ is trivial 
since $A\in M_{n\times k}(\Z )$.)  

On the other hand, if necessary, by replacing 
a basis of $\Z^n$, $A$ can be of the form 
\[
A=
\bordermatrix{
       &        &  k                       &          \cr
       & d_1 &                           &         \cr 
\ \ k  &        & \ddots               &         \cr
       &        &                            & d_k \cr 
n-k  &        & \boldsymbol{0} &           
} , 
\]  
and
\begin{align}
{\rm Im}\varphi_A& =\overbrace{\R\oplus \cdots \oplus \R}^k 
\oplus \{ 0\} \oplus \cdots \oplus \{ 0\}
\subset \R^n
\nonumber \\
{\rm Im}\varphi_A\cap \Z^n&=\overbrace{\Z\oplus \cdots \oplus \Z}^k 
\oplus \{ 0\} \oplus \cdots \oplus \{ 0\}
\subset \Z^n \nonumber \\
\varphi_A(\Z^k)&=d_1\Z\oplus \cdots \oplus d_k\Z 
\oplus \{ 0\} \oplus \cdots \oplus \{ 0\} \subset \Z^n . \nonumber 
\end{align}
(Use the fundamental divisor theory. )  Then 
\begin{align}
\therefore \ \ 
{\rm Im}\varphi_A\cap \Z^n \subset \varphi_A(\Z^k) 
& \Longleftrightarrow d_i=\pm 1\ (i=1, \ldots ,k) \nonumber \\
& \Longleftrightarrow \{ u_1, \ldots , u_k\} {\rm : primitive} . 
\nonumber 
\end{align}
\end{proof}

Assume $\Phi^{-1}(\lb )\neq \emptyset$. For each 
$(\xi , \theta , z )\in \Phi^{-1}(\lb )$, we set 
\[
I_{(\xi , \theta , z )}= \{ i\in \{ 1, \ldots ,d\} \lvert z_i=0 \} . 
\]
\begin{lem}\label{a2}
$(1)$ The following conditions are equivalent. 
\begin{enumerate}
\item $\lambda$ is a regular value of $\Phi$.   
\item For arbitrary $(\xi , \theta , z )\in \Phi^{-1}(\lb )$, 
$\{ u_i\}_{i\in I_{(\xi , \theta , z )}}$ is linear independent. 
\end{enumerate}
$(2)$ Under the condition of $(1)$, the following conditions are equivalent. 
\begin{enumerate}
\item $T^k$-action on $\Phi^{-1}(0)$ is free.   
\item For arbitrary $(\xi , \theta , z )\in \Phi^{-1}(\lb )$, 
$\{ u_i\}_{i\in I_{(\xi , \theta , z )}}$ is primitive. 
\end{enumerate}
\end{lem}
\begin{proof}
$\Phi$ can be decomposed into two maps $\mu_{\C^d}:\C^d \to \R^d$ in Example \ref{1.2} 
and $\Phi_1:T^*T^n\to \R^d$ which is defined by 
\[
\Phi_1(\xi , \theta )=\ ^tB\xi , 
\]
where $B=(u_1, \ldots , u_d)\in M_{n\times d}(\Z)$. 
Let us pay attention to $\mu_{\C^d}$. 
Since $\mu_{\C^d}$ is the moment map of $T^d$-action on $\C^d$, we have 
${\rm Im}(d(\mu_{\C^d})_z)^{\perp}={\rm Lie}(T^d_{z})$, 
where $T^d_z$ is the stabilizer of $z$ of $T^d$-action on $\C^d$.   
Since $T^d_z$ can be written by 
\[
T^d_{z}=\{ t \in T^d \colon 
i\notin I_{(\xi , \theta , z )} \Rightarrow t_i=0 \} , 
\]
we have 
\begin{equation}\label{eq1.2}
{\rm Im}(d(\mu_{\C^d})_z)^{\perp}=\{ (\lambda_1, \ldots ,\lambda_d)\in \R^d 
\colon i\notin I_{(\xi , \theta , z )} \Rightarrow \lambda_i=0\} . 
\end{equation}
On the other hand, $d(\Phi_1)_{(\xi , \theta )}:T_{(\xi , \theta )}(T^*T^n)\cong 
T_{\xi}\R^n \oplus T_{\theta}T^n \to \R^d$ is 
\begin{equation}\label{eq1.3}
d(\Phi_1)_{(\xi , \theta )}(v_{\R^n}, v_{T^n})=\ ^tBv_{\R^n} . 
\end{equation}
From (\ref{eq1.2}) and (\ref{eq1.3}), we have 
\[
\text{$(\xi , \theta , z)$ is a regular point of $\Phi$}
\Leftrightarrow \{ u_i\}_{i\in I_{(\xi , \theta , z)}} 
\text{are linear independent}. 
\]
This proves (1). 

(2) is obtained from the fact that the stabilizer $T^d_{(\xi , \theta ,z)}$ of 
$(\xi ,\theta ,z)$ of $T^d$-action on $T^*T^n\times \C^d$ is written by
\[
T^d_{(\xi , \theta , z)}=\{ t \in T^d \colon  
\sum_{i\in I_{(\xi , \theta , z)}}t_iu_i\in \Z^n \} 
\] 
and Lemma \ref{a1}. 
\end{proof}
\begin{thm}\label{a3}
$\overline{(T^*T^n)}_{\{ \mu_{u_i}\ge \lambda_i \}_{i=1, \ldots ,d}}$ 
is smooth, if and only if for each $[\xi , \theta , z]\in 
\overline{(T^*T^n)}_{\{ \mu_{u_i}\ge \lambda_i \}_{i=1, \ldots ,d}}$, 
$\{ u_i\}_{i\in I_{(\xi, \theta , z)}}$ is primitive. 
\end{thm}
\begin{proof}
This is a direct consequence of Lemma \ref{a2}. 
\end{proof}

By Remark \ref{1.10}, $\overline{(T^*T^n)}_{\{ \mu_{u_i}\ge \lambda_i \}_{i=1, \ldots ,d}}$ 
is equipped with the Hamiltonian $T^n$-action which is induced from the natural 
$T^n$-action on $T^*T^n$. 
\begin{thm}\label{a4}
Under the condition in Theorem $\ref{a3}$, in the particular case 
where $d\le n$ and $\lambda_1=\cdots =\lambda_d=0$, there exists an automorphism 
$\rho$ of $T^n$ such that $\overline{(T^*T^n)}_{\{ \mu_{u_i}\ge 0 \}_{i=1, \ldots ,d}}$ is 
$\rho$-equivariantly symplectomorphic to the $T^n$-action on 
$(\C^d\times T^*T^{n-d},\om_{\C^d}\oplus \om_{T^*T^{n-d}})$ 
that is the direct product of the $T^d$-action on $\C^d$ in 
Example $\ref{1.2}$ and the $T^{n-d}$-action on $T^*T^{n-d}$ in Example $\ref{1.3}$. 
\end{thm}
\begin{proof}
For primitive 
$u_1, \ldots , u_d$, there exists an element $\rho \in \GL_n(\Z )$ such that $\rho$ 
maps $u_i$ to the $i$th fundamental vector $e_i$ for $i=1, \ldots ,d$. Then, $\rho$ 
induces the symplectomorphism $\widetilde{\varphi}=\ ^t\rho^{-1}\times \rho$ of 
$(T^*T^n=\R^n\times T^n, \om_{T^*T^n})$. The source $T^*T^n$ is equipped with 
$d$ commutative Hamiltonian circle actions defined by $u_1$, $\ldots$, $u_d$ with 
moment maps $\mu_{u_1}$, $\ldots$, $\mu_{u_d}$, and on the target $T^*T^n$, 
$e_1$, $\ldots$, $e_d$ define $d$ commutative Hamiltonian circle actions with 
moment maps $\mu_{e_1}$, $\ldots$, $\mu_{e_d}$, respectively. Since for each $i$, 
$\widetilde{\varphi}$ is equivariant with respect to the Hamiltonian circle action defined 
by $u_i$ on the source and that defined by $e_i$ on the target, $\widetilde{\varphi}$ 
descends to the $\rho$-equivariantly symplectomorphism $\varphi$ from the simultaneous 
cut space $\overline{(T^*T^n)}_{\{ \mu_{u_i}\ge 0 \}_{i=1, \ldots ,d}}$ to 
$\overline{(T^*T^n)}_{\{ \mu_{e_i}\ge 0 \}_{i=1, \ldots ,d}}$. But by Example 
\ref{1.11},  $\overline{(T^*T^n)}_{\{ \mu_{e_i}\ge 0 \}_{i=1, \ldots ,d}}$ is 
equivariantly symplectomorphic to $\C^d\times T^*T^{n-d}$. This proves Theorem \ref{a4}. 
\end{proof}

\section{Smoothness of induced maps between cut spaces}
Let $u, v\in \Z^n$ be two vectors each of which is primitive, $\rho \in SL_n (\Z )$ 
with $\rho (u)=\pm v$, and $\varphi :U\to V$ an orientation preserving diffeomorphism 
between two open sets $U$, $V\subset \R^n$which satisfies 
\[
\varphi \left(\{ \xi \in U\colon \< u, \xi \> \ge 0\} \right) 
=\{ \eta \in V\colon \< v, \eta \> \ge 0\} . 
\]  
In this Appendix, we shall show the following lemma. 
\begin{lem}\label{B1}
If $\varphi$ satisfies the condition 
\[
\< u, \xi \> =\< v, \varphi (\xi )\>
\]
on a sufficiently small neighborhood of $\{ \xi \in U \colon \< u, \xi \> =0\}$ in 
$\{ \xi \in U\colon \< u, \xi \> \ge 0\}$, the map 
$\varphi \times \rho :U\times T^n\to V\times T^n$ descends to an orientation preserving 
diffeomorphism between cut spaces $\overline{(U\times T^n)}_{\mu_u\ge 0}$ and 
$\overline{(V\times T^n)}_{\mu_v\ge 0}$. 
\end{lem}
\begin{proof}
Let $\Phi_u:U\times T^n\times \C \to \R$ and $\Phi_v:V\times T^n\times \C \to \R$ 
be the maps which are defined by 
\[
\Phi_u(\xi ,\theta ,z)=\< u,\xi \>-\lvert z\rvert^2, \ \ 
\Phi_v(\eta ,\tau ,w)=\< v,\eta \>-\lvert w\rvert^2. 
\]
Define the diffeomorphism $\overline{\psi}:\Phi_u^{-1}(0)\to \Phi_v^{-1}(0)$ by 
\[
\overline{\psi}(\xi ,\theta ,z)=\left( \varphi (\xi ), \rho (\theta ), 
\sqrt{\frac{\< v,\varphi (\xi )\>}{\< u,\xi \>}} \overline{z}^{\rho}\right) , 
\]
where 
\[
\overline{z}^{\rho}=
\begin{cases}
z & \text{if } \rho (u)=v \\
\overline{z} & \text{if } \rho(u)=-v .
\end{cases}
\]
Note that $\overline{\psi}$ is well-defined from the assumption of the lemma. 
$\overline{\psi}$ is equivariant with respect to the free circle actions 
\[
t\cdot (\xi ,\theta ,z)=(\xi ,\theta +tu, e^{-2\pi\sqrt{-1}t}z)
\] 
on $\Phi_u^{-1}(0)$ and 
\[
t\cdot (\eta ,\tau ,w)=(\eta ,\tau +t\rho (u), \overline{e^{-2\pi\sqrt{-1}t}}^{\rho}w) 
\]
on $\Phi_v^{-1}(0)$. Then $\overline{\psi}$ descends to the orientation preserving 
diffeomorphism $\psi : \Phi_u^{-1}(0)/S^1\to \Phi_v^{-1}(0)/S^1$ since 
$\varphi$ and $\rho$ preserve orientations. 
\end{proof}
By using Lemma \ref{B1} repeatedly, we can show the version of the lemma for 
the case of simultaneous cut spaces. Let $\{u_1, \ldots , u_d\}$, 
$\{ v_1, \ldots ,v_d\} \subset \Z^n$ be the primitive tuples of vectors, 
$\rho \in SL_n(\Z )$ with $\rho (u_i)=\pm v_i$ for $i=1,\ldots , d$, and 
$\varphi :U\to V$ an orientation preserving diffeomorphism 
between two open sets $U$, $V\subset \R^n$which satisfies 
\[
\varphi \left(\{ \xi \in U\colon \< u_i, \xi \> \ge 0\ i=1,\ldots ,d \} \right) 
=\{ \eta \in V\colon \< v_i, \eta \> \ge 0\ i=1,\ldots ,d\} . 
\]  
\begin{lem}\label{B2}
If $\varphi$ satisfies the condition 
\[
\< u_i, \xi \> =\< v_i, \varphi (\xi )\>
\]
on a sufficiently small neighborhood of 
$\{ \xi \in U \colon \< u_i, \xi \> =0,\ \< u_j, \xi \> \ge 0\ j\neq i \}$ in 
$\{ \xi \in U\colon \< u_j, \xi \> \ge 0\ j=1,\ldots ,d \}$ for each $i=1,\ldots ,d$, 
the map $\varphi \times \rho :U\times T^n\to V\times T^n$ descends to the 
orientation preserving diffeomorphism between simultaneous cut spaces 
$\overline{(U\times T^n)}_{\{ \mu_{u_i}\ge 0\}_{i=1,\ldots ,d}}$ and 
$\overline{(V\times T^n)}_{\{ \mu_{v_i}\ge 0\}_{i=1,\ldots ,d}}$. 
\end{lem}

\vspace*{5mm}
Graduate School of Mathematical Sciences, The University of Tokyo, \\
8-1 Komaba 3-chome, Meguro-ku, Tokyo, 153-8914, Japan \\
e-mail:takahiko@ms.u-tokyo.ac.jp

\end{document}